\documentclass[11pt]{amsart}
\usepackage{amsmath,amssymb,amsthm}
\usepackage{float}
\usepackage{tikz-cd}
\usepackage[shortlabels]{enumitem}

\newtheorem{theorem}{Theorem}
\newtheorem*{theorem*}{Theorem}
\newtheorem{proposition}[theorem]{Proposition}
\newtheorem{notation}[theorem]{Notation}
\newtheorem{lemma}[theorem]{Lemma}

\newtheorem{corollary}[theorem]{Corollary}

\numberwithin{theorem}{section}
\numberwithin{equation}{section}

\theoremstyle{definition}
\newtheorem{definition}[theorem]{Definition}
\newtheorem{remark}[theorem]{Remark}
\newtheorem{example}[theorem]{Example}

\makeatletter
\def\thmhead@plain#1#2#3{%
  \thmname{#1}\thmnumber{\@ifnotempty{#1}{ }\@upn{#2}}%
  \thmnote{ {\the\thm@notefont#3}}}
\let\thmhead\thmhead@plain
\makeatother

\def\C{\mathbb{C}}

\def\R{\mathbb{R}}
\def\Z{\mathbb{Z}}

\def\sC{\mathcal{C}}
\def\sD{\mathcal{D}}
\def\sE{\mathcal{E}}

\def\sH{\mathcal{H}}

\def\sL{\mathcal{L}}

\def\sR{\mathcal{R}}

\def\sW{\mathcal{W}}

\def\fB{\mathfrak{B}}

\def\fX{\mathfrak{X}}

\def\bra{\langle}
\def\ket{\rangle}

\newcommand{\spacedt}[1]
{
\quad{\textup{#1}}\quad
}

\newcommand{\Matrix}[1]
{
\left(
\begin{matrix}
#1
\end{matrix}
\right)
}

\newcommand{\restrto}[2]
{
\left.{#1}\right|_{#2}
}

\newcommand{\su}[1]
{
_{_{#1}}
}

\newcommand{\norm}[1]
{
\left\|#1\right\|
}

\def\RA{\Rightarrow}
\def\imm{\hookrightarrow}

\def\Op{\mathcal{O}p}
\def\abadapted{$(\alpha,\beta)$-adapted}

\title{Riemannian properties of Engel structures}

\author{Nicola Pia}
\date{The author is supported by the DAAD programme Research Grants for Doctoral Candidates and Young Academics and Scientists (more than 6 months) No. 57381410, 2018/19}
\address{Mathematisches Institut, LMU M\"unchen, Theresienstr. 39, 80333 M\"unchen, Germany}
\email{pia@math.lmu.de}

\begin{document}
\begin{abstract}
This paper is about geometric and Riemannian properties of Engel structures, i.e. maximally non-integrable $2$-plane fields on $4$-manifolds.
Two $1$-forms $\alpha$ and $\beta$ are called Engel defining forms if $\mathcal{D}=\ker\alpha\cap\ker\beta$ is an Engel structure and $\mathcal{E}=\ker\alpha$ is its associated even contact structure, i.e. $\mathcal{E}=[\mathcal{D},\mathcal{D}]$.
A choice of Engel defining forms determines a distribution $\mathcal{R}$ transverse to $\mathcal{D}$ called the Reeb distribution.
We study conditions that ensure integrability of $\mathcal{R}$.
For example if we have a metric $g$ which makes the splitting $TM=\mathcal{D}\oplus\mathcal{R}$ orthogonal and such that $\mathcal{D}$ is totally geodesic then there exists an integrable Reeb distribution $\tilde{\mathcal{R}}$.

It turns out that integrabilty of $\mathcal{R}$ is related to the existence of vector fields $Z$ whose flow preserves $\mathcal{D}$, so called Engel vector fields.
A K-Engel structure is a triple $(\mathcal{D},\,g,\,Z)$ where $\mathcal{D}$ is an Engel structure, $g$ is a Riemannian metric, and $Z$ is a vector field which is Engel, Killing, and orthogonal to $\mathcal{E}$.
In this case we can construct Engel defining forms with very nice properties and such that $\mathcal{R}$ is integrable.
Moreover we can classify the topology of K-Engel manifolds studying the action of the flow of $Z$.
As natural consequences of these methods we provide a construction which is the analogue of the Boothby-Wang construction in the contact setting and we give a notion of contact filling for an Engel structure.
\end{abstract}
\maketitle

\section{Introduction}
The theory of contact metric structures, K-contact and Sasakian manifold is an active field of research (see \cite{blair, boyer} for an introduction).
Moreover the study of the interplay between metrics and contact structures has helped in the understanding of topological properties of these structures \cite{ekm, massot}.
The goal of this work is to study analogous problems for Engel structures.

An Engel structure on a smooth $4$-manifold is a smooth maximally non-integrable $2$-plane distribution $\sD$, i.e. $\sE:=[\sD,\sD]$ is a rank $3$ distribution such that $[\sE,\sE]=TM$.
The distribution $\sE$ is an even contact structure, which is the even dimensional analogue of a contact structure.
If for a given Engel structure $\sD$ we have two $1$-forms $\alpha$ and $\beta$ satisfying $\sD=\ker\alpha\cap\ker\beta$ and $\sE=\ker\alpha$, we say that $\alpha$ and $\beta$ are Engel defining forms.
This happens if and only if
\[
	\alpha\wedge d\alpha\ne0,\qquad \alpha\wedge\beta\wedge d\beta\ne0\spacedt{and}\alpha\wedge d\alpha\wedge\beta=0,
\]
A pair of defining forms determines a distribution $\sR=\bra T,\, R\ket$ transverse to $\sD$ via
\[
\begin{matrix}
  i_T(\alpha\wedge d\beta)=0,&\beta(T)=1,&\alpha(T)=0, \\
  i_R(\beta\wedge d\beta)=0,&\beta(R)=0,&\alpha(R)=1.
\end{matrix}
\]
This is called the \emph{Reeb distribution associated with} $\alpha$ and $\beta$.

Engel defining forms first appeared in \cite{adachi}, but they were already used as a technical tool in \cite{montgomeryPaper}.
Moreover the vector fields $T$ and $R$ were introduced in \cite{adachiWeinstein}, but the properties of the distribution $\sR$ were never studied explicitly.
One of our goals is to understand geometric properties of $\alpha$ and $\beta$ and of their associated Reeb distribution.

A natural question to ask is whether we can choose Engel defining forms so that the associated Reeb distribution is integrable.
We present many ways of rephrasing this problem one of which gives it a geometrical flavour.
The \emph{contactization} of the Engel structure $(M,\, \sD=\ker\alpha\wedge\beta)$ is the contact manifold $(X=M\times\R,\, \eta=\beta+s\alpha)$.
The existence of an integrable Reeb distribution is linked with how the Reeb vector field of $(X,\,\eta)$ intersects graphical hypersurfaces.

The solution to this problem in general is not known, but it turns out that some metric properties ensure a positive answer in many cases.
Since a choice of Engel defining forms $\alpha$ and $\beta$ determines uniquely the splitting $TM=\sD\oplus\sR$, it is natural to consider Riemannian metrics $g$ which satisfy $\sD^\bot=\sR$.
\begin{theorem*}
	Let $\sD=\ker\alpha\wedge\beta$ be an Engel structure with a choice of Engel defining forms, and let $g$ be a metric such that $\sD^\bot=\sR$.
	If $\sD$ is totally geodesic, then there exists a nowhere vanishing function $\mu\in\sC^\infty(M)$ such that $\tilde\beta=\mu\beta$ satisfies $d\tilde\beta^2=0$.
	In particular, $\tilde\sR$ is integrable.
\end{theorem*}

Other conditions for integrability are given by the existence of \emph{Engel vector fields} (i.e.  vector fields whose flow preserves $\sD$) with some additional properties.
A \emph{K-Engel structure} is an Engel structure $(M,\sD)$ together with a metric $g$ and a vector field $Z$ which is Engel, Killing and orthogonal to $\sE$.
The existence of such a vector field allows the construction of defining forms $\alpha$ and $\beta$ satisfying $d\alpha^2=0=d\beta^2$ and $\beta\wedge d\alpha=0$.
This in turn provides a framing $\fB=\{W,\, X,\, T,\, R\}$ for $TM$ with $\sW=\bra W\ket$, $\sD=\bra W,\, X\ket$ and such that $R=Z$ commutes with all other vector fields in the framing.
The closure of the flow of $R$ in the group of isometries of $(M,g)$ is a torus $T^k$ acting on $M$.
The existence of the framing $\fB$ ensures that this torus acts freely, so that we have the following result.
\begin{theorem*}
	If $(M,\sD)$ admits a K-Engel structure, then $M$ is diffeomorphic to one of the following:
	\begin{itemize}
		\item $T^4$;
		\item a principal $T^2$-bundle over a surface;
		\item a principal $S^1$-bundle over a $3$-manifold.
	\end{itemize} 
\end{theorem*}
The dimension $k$ of the closure of the orbits of $R$ is the rank of the K-Engel structure.

The most interesting case is rank $1$, because all of these examples come from the Engel Boothby-Wang construction.
This is the analogue of the Boothby-Wang construction of K-contact manifolds, and it already appeared, in a different context, in \cite{mitsu} under the name of prequantum prolongation.
A well-known theorem in K-contact geometry asserts that every K-contact structure can be perturbed to a quasi-regular one, the following result is the analogue for K-Engel geometry.
\begin{theorem*}
 Let $(M,\,\sD=\ker\alpha\wedge\beta,\,g,\,R)$ be a K-Engel structure of rank $k$, then there exist K-Engel structures of rank $1$ $(M,\,\sD=\ker\alpha_i\wedge\beta_i,\,g_i,\,R_i)$ for $i=1,...,k$ such that $R_1,...,R_k$ are linearly independent.
\end{theorem*}

Another remarkable property of rank $1$ structures is that they all admit a \emph{contact filling}.
This means that they can be realized as boundaries of a contact manifold $(X,\,\eta)$ so that a collared neighbourhood of $\partial X$ is isomorphic to the \emph{contactization} of $M$.
This is the analogue of a symplectic filling for a contact structures.

\subsection{Structure of the Paper}
In Section~\ref{SEC_EvenContact} we recall some useful facts from the theory of even contact structures and in Section~\ref{SEC_EngelForms} we introduce Engel defining forms $\alpha$ and $\beta$.
The end of Section~\ref{SEC_EngelForms} as well as Section~\ref{SEC_IntegrableR} concerns conditions that ensure the existence of defining forms such that the associated Reeb distribution is integrable.
Some of these conditions do not give any geometric insight, but are just reformulations which are useful for the remainder of the chapter.
We present a condition that is geometric, and that is linked to the contactization of an Engel structure.

In Section~\ref{SEC_MetricProperties} we study when $\sD=\ker\alpha\wedge\beta$ and its associated Reeb distribution $\sR$ are totally geodesic.
We prove that if $g$ is a metric such that $\sD^\bot=\sR$ then $\sR$ cannot be totally geodesic.
Moreover if $\sD$ is totally geodesic then there exists a multiple $\tilde\beta$ of $\beta$ such that $d\tilde\beta^2=0$.
This is a sufficient condition for the existence of a Reeb distribution which is integrable.

In Section~\ref{SEC_EngelVF} we recall the theory of \emph{Engel vector fields} $Z$, i.e. vector fields satisfying $\sL_Z\sD\subset\sD$.
Many of the results we present were developed in \cite{vogelThesis}, but we rephrase them in terms of the framing induced by a pair of defining forms $\alpha$ and $\beta$.
In Section~\ref{SEC_EngelKillingVF} we specialize to the case where $Z$ is also Killing with respect to a metric $g$ that makes it orthogonal to $\sE$.
Under these hypothesis we can find defining forms verifying $d\alpha^2=0=d\beta^2$ and $d\alpha\wedge\beta=0$.
This allows the construction of a framing $\{W,\, X,\, T,\, R\}$ such that $R$ commutes with every other vector field in the framing.

In Sections~\ref{SEC_rho} and~\ref{SEC_TopologyOfKEngel} we analyse the topology of manifolds that admit a K-Engel structure.
We study more in details the bundle structure of the manifolds that appear in this classification, and we exhibit explicit constructions of K-Engel structures.
The most remarkable one is the \emph{Engel Boothby-Wang} construction, the analogue of the one in the K-contact case.

In Section~\ref{SEC_ContactFillings} we define \emph{contact fillings} of an Engel structure.
The construction is inspired by symplectic fillings of a contact structure, but we do not know if it implies any rigidity as in the contact case.
We prove that Engel Boothby-Wang manifolds are examples of Engel structures which admit contact fillings.

Finally in Section~\ref{SEC_GeoEngel} we go through the list of geometric Engel structures (see \cite{vogelGeometricEngel}) and we verify that almost all of them admit a left invariant K-Engel structure.
This provides many examples of K-Engel manifolds.

\textbf{Acknowledgements:}
I would like to thank my advisors Prof Gianluca Bande for being always there for me and motivating me throughout my PhD studies; and Prof Dieter Kotschick for the innumerable discussions and inspiring questions.
I would also like to thank Prof Thomas Vogel and Prof Vincent Colin for their useful feedback.
Finally I need to thank Rui Coelho and Giovanni Placini for having the patience to listen to me all the time.

\section{Some properties of even contact structures}\label{SEC_EvenContact}
In what follows all manifolds are supposed to be smooth and orientable.

An \emph{even contact structure} $\sE$ is a maximally non-integrable hyperplane field on an even dimensional manifold $M^{2n+2}$.
This means that $\sE$ is locally the kernel of a $1$-form $\alpha$ satisfying $\alpha\wedge d\alpha^{2n}\ne0$.
An even contact structure is co-orientable if  orientable, this is equivalent to the existence of a (global) $1$-form $\alpha$ such that $\sE=\ker\alpha$.
There exists a unique line field $\sW\subset\sE$ whose holonomy preserves $\sE$, i.e. $[\sW,\sE]=\sE$, we call $\sW$ \emph{characteristic, kernel or Cauchy line field of} $\sE$.
If $\sE=\ker\alpha$ then $\sW=\ker\alpha\wedge d\alpha^{2n}$.

We will be interested in properties of vector fields whose flow preserves $\sE$.
The existence of such vector fields depends on the dynamics of the characteristic foliation.
The following result (which appeared in \cite{biengel}) ensures that the holonomy of $\sW$ is volume-preserving if and only if there is a transverse symmetry for $\sE$.
\begin{proposition}[\cite{biengel}]\label{PROP_PropertiesOfVolumePreservingW}
	Let $(M^{2n+2},\,\sE=\ker\alpha)$ be an even contact structure with orientable characteristic foliation $\sW=\ker\alpha\wedge d\alpha^n$, then the following are equivalent:
	\begin{enumerate}[1.]
		\item $\sW$ has volume-preserving holonomy;
		\item $\sW$ is the kernel of a closed $(2n+1)$-form;
		\item $\alpha$ can be chosen so that $d\alpha$ has constant rank $2n$;
		\item there exists a vector field $Z\in\fX(M)$ transverse to $\sE$ whose flow preserves $\sE$.
	\end{enumerate}
\begin{proof}
	The equivalence between $1$, $2$, and $3$ was proved in \cite{biengel}, for more details see also \cite{me}.
	\fbox{$3\RA4$} If $d\alpha$ has constant rank $2n$ then $\ker d\alpha$ is a $2$-plane field.
	Now $W$ is in the kernel of $d\alpha$ and $\ker d\alpha$ is trivial as a bundle.
	Let $\ker d\alpha=\bra W,\, Z\ket$.
	Uniqueness of $\sW$ ensures that $Z$ must be transverse to $\sE$.
	Choose $Z\in\ker d\alpha$ such that $\alpha(Z)=1$, then
	\[
		\sL_Z\alpha=i_Zd\alpha+di_Z\alpha=d(\alpha(Z))=0.
	\]
	
	\fbox{$4\RA3$} By hypothesis, for any defining form $\alpha$ we have $\sL_Z\alpha=\lambda\alpha$.
	Pick $\alpha$ so that $\alpha(Z)=1$. 
	Taking the derivative of this we get
	\[ 
	  0=\sL_Z(\alpha(Z))=(\sL_Z\alpha)(Z)+\alpha(\sL_ZZ)=\lambda,
	\]
	hence $i_Zd\alpha=\sL_Z\alpha=0$, and the rank of $d\alpha$ is never maximal, so it must be $2n$ everywhere.
\end{proof}
\end{proposition}

\begin{remark}
	If the flow of a vector field $Z\in\fX(M)$ preserves the even contact structure $\sE$, then it must also preserve its characteristic foliation $\sW$.
	If we are in the hypothesis of Proposition~\ref{PROP_PropertiesOfVolumePreservingW} this implies in particular that, if $W$ is a section of $\sW$, we have $\sL_ZW=aW$ and $\bra W,Z\ket$ is a foliation.
	Another way of seeing this is taking $\alpha$ defining form such that $\alpha(Z)=1$, since this ensures $\ker d\alpha=\bra W,Z\ket$.
\end{remark}

\section{Engel defining forms and the Reeb distribution}\label{SEC_EngelForms}

An \emph{Engel structure} $\sD$ is a smooth $2$-plane field on a smooth $4$-manifold such that $\sE:=[\sD,\sD]$ is an even contact structure.
One can proof that the characteristic foliation $\sW$ of $\sE$ satisfies $\sW\subset\sD$ hence giving the flag of distributions
\[
 \sW\subset\sD\subset\sE\subset TM,
\]
called \emph{Engel flag}.
If all distributions in the flag are orientable then $TM$ is trivial.
This is a strong contraint on the topology of $M$ (see \cite{montgomeryPaper} for more details).
The question of existence of Engel structures on parallelizable $4$-manifolds has been studied in \cite{cp3, overtwistedEngel, vogelExistence}.
We will instead focus the attention on the study of Engel defining forms.

\begin{definition}\label{DEF_EngelForms}
	Let $\sD$ be an Engel structure.
	Two $1$-forms $\alpha,\,\beta$ are called \emph{Engel defining forms} if they verify $\sD=\ker\alpha\cap\ker\beta$ and $\sE=\ker\alpha$.
	Equivalently
	\begin{subequations}
		\label{EQN_DEF}
		\begin{align}
		\alpha\wedge d\alpha\ne0\label{EQN_DEF1} \\
		\alpha\wedge\beta\wedge d\beta\ne0 \label{EQN_DEF2} \\ 
		\alpha\wedge d\alpha\wedge\beta=0 \label{EQN_DEF3}
		\end{align}
	\end{subequations}
\end{definition}

Observe that Equation~\eqref{EQN_DEF1} ensures that $\sE=\ker\alpha$ is an even contact structure; denote with $\sW=\ker\alpha\wedge d\alpha$ its characteristic foliation.
Equation~\eqref{EQN_DEF2} implies that $\ker\beta$ is an even contact structure, and it ensures that its characteristic foliation is transverse to $\sE$.
Finally Equation~\eqref{EQN_DEF3} implies that $\sW\subset\ker\beta$.

\begin{remark}\label{REM_EngelasIntersOfEvenContact}
	An orientable Engel structure $\sD$ on an orientable manifold $M$ can always be seen as the intersection of two even contact structures $\sD=\sE\cap\sE'$.
	The first one is uniquely determined by the condition $[\sD,\sD]=\sE$, and the second one instead is not unique, and must satisfy $\sW\subset\sE'$ and $\sW'\pitchfork\sE$.
	A choice of Engel defining forms $\alpha$ and $\beta$ corresponds to  choice of a co-orientation $\sE=\ker\alpha$ and to a choice of $\sE'$ together with a co-orientation $\sE'=\ker\beta$.
\end{remark}

\begin{lemma}\label{LEM_ExistenceOfDefiningForms}
	Let $M^4$ be parallelizable then $\sD$ Engel structure on $M$ admits Engel defining forms if and only if $\sD$ is orientable.
	\begin{proof}
		If $\alpha$ and $\beta$ are Engel defining forms then $d\beta$ is an orientation of $\sD$.
		Conversely, suppose that $\sD$ is orientable.
		Since $M$ is orientable then $\sW$ is trivial as a line bundle, let $W$ be a non-singular section.
		Since $\sD$ is orientable, it has a non-singular section $X$ nowhere tangent to $\sW$.
		Now $Y=[W,X]$ is a non-singular section of $\sE$ nowhere tangent to $\sD$, and $Z=[X,Y]$ is transverse to $\sE$.
		Choose $\alpha$ such that $\ker\alpha=\bra W,\, X,\, Y\ket$ and $\alpha(Z)=1$, and similarly $\beta$ such that $\ker\beta=\bra W,\, X,\, Z\ket$ and $\beta(Y)=1$.
	\end{proof}
\end{lemma}

Notice that if $\alpha$ and $\beta$ are defining forms for $\sD$, then all other possible defining forms are given by
\[
\tilde\alpha=\lambda\alpha,\qquad\tilde\beta=\mu\beta+\nu\alpha
\]
where $\lambda,\mu,\nu\in\sC^\infty(M)$ and $\lambda,\,\mu$ are nowhere vanishing.

In what follows we will write $\sD=\ker\alpha\wedge\beta$ to denote the choice of Engel defining forms $\alpha$ and $\beta$ for the Engel structure $\sD$.

\subsection{Reeb distribution}
A choice of Engel defining forms determines uniquely a distribution $\sR$ transverse to $\sD$.
Indeed Equation~\eqref{EQN_DEF2} implies that $\alpha\wedge d\beta$ and $\beta\wedge d\beta$ are nowhere vanishing $3$-forms, in particular they have a $1$-dimensional kernel.
Take $T$ nowhere vanishing section of $\ker\alpha\wedge d\beta$, this implies $\alpha(T)=0$.
On the other hand Equation~\eqref{EQN_DEF2} ensures that $\beta(T)\ne0$, hence we can normalize $T$ via $\beta(T)=1$.
Similarly pick $R$ nowhere vanishing section of $\ker\beta\wedge d\beta$ and normalize it via $\alpha(R)=1$.
By construction we have $\beta(R)=0$.
\begin{definition}\label{EQN_ReebDistribution}
	Let $(M,\sD=\ker\alpha\wedge\beta)$ be an Engel structure with a choice of defining forms $\alpha$ and $\beta$.
	The \emph{Reeb distribution} associated with $\alpha$ and $\beta$ is $\sR=\bra T,\,R\ket$ where
	\[
	\begin{matrix}
	i_T(\alpha\wedge d\beta)=0,&\beta(T)=1,&\alpha(T)=0,\\
	i_R(\beta\wedge d\beta)=0,&\beta(R)=0,&\alpha(R)=1. 
	\end{matrix}
	\]
\end{definition}

The definition implies $TM=\sD\oplus\sR$.
We are interested in understanding geometric properties of $\sR$.
In general the Reeb distribution is not integrable.
The following result gives a necessary and sufficient condition for integrability.
\begin{proposition}\label{PROP_IntegrabilityOfR}
	Let $\sD=\ker\alpha\wedge\beta$ be an Engel structure and consider the associated Reeb distribution $\sR=\bra T,\,R\ket$, and denote by $c\su{TR}=d\beta(R,T)$.
	Then $\sR=\ker(d\beta+c\su{TR}\beta\wedge\alpha)$ and it is integrable if and only if
	\begin{equation}\label{EQN_integrabilityCondR}
	d(c\su{TR}\alpha)\wedge\beta=0.
	\end{equation}
	\begin{proof}
		For the first assertion we notice that $i_R(\beta\wedge d\beta)=0$ and $\beta(R)=0$ imply that $i_Rd\beta$ is a multiple of $\beta$.
		In fact $i_Rd\beta=c\su{TR}\beta$, so that $R\in\ker(d\beta+c\su{TR}\beta\wedge\alpha)$.
		Similarly we have have $T\in\ker(d\beta+c\su{TR}\beta\wedge\alpha)$ so that $\sR\subset\ker(d\beta+c\su{TR}\beta\wedge\alpha)$.
		On the other hand the kernel of this $2$-form has rank $2$ because
		\[
		\alpha\wedge\beta\wedge(d\beta+c\su{TR}\beta\wedge\alpha)=\alpha\wedge\beta\wedge d\beta\ne0.
		\]
		
		To prove the last claim choose two $1$-forms $\rho,\tau$ such that
		\begin{equation}\label{EQN_effimeral}
		\rho\wedge\tau=d\beta+c\su{TR}\beta\wedge\alpha,
		\end{equation}
		which is possible since $\sD$ is orientable.
		By Frobenius' Theorem $\sR$ is integrable if and only if $\rho\wedge\tau\wedge d\tau=0$ and $\tau\wedge\rho\wedge d\rho=0$.
		Differentiating~\eqref{EQN_effimeral} we get
		\[
		d\rho\wedge\tau-\rho\wedge d\tau=c\su{TR}d\beta\wedge\alpha -d(c\su{TR}\alpha)\wedge\beta=c\su{TR}\rho\wedge\tau\wedge\alpha-d(c\su{TR}\alpha)\wedge\beta,
		\]
		hence the integrability condition translates to 
		\begin{equation}\label{daf}
		d(c\su{TR}\alpha)\wedge\beta\wedge\tau=0\spacedt{and}d(c\su{TR}\alpha)\wedge\beta\wedge\rho=0.
		\end{equation}
		Since $d(c\su{TR}\alpha)\wedge\beta\wedge\alpha=c\su{TR}d\alpha\wedge\beta\wedge\alpha=0$ and $d(c\su{TR}\alpha)\wedge\beta\wedge\beta=0$, conditions~\eqref{daf} are satisfied if and only if $d(c\su{TR}\alpha)\wedge\beta=0$.
	\end{proof}
\end{proposition}

\subsection{A useful technical lemma}

In this section we list some formulas which describe Lie brackets of a framing induced by $\sD=\ker\alpha\wedge\beta$.

\begin{notation}\label{NOT_Commutators}
	In what follows we will fix a parallelization $\{W,\, X,\, T,\, R\}$ of $M$ and we will use the letters $a,\, b,\, c$ and $d$ to denote respectively the $W,\, X,\, T$ and $R$ components of vector fields.
	Moreover we will use lower indices $a_{AB}, b_{AB}, c_{AB}$ and $d_{AB}$ to denote the components of the Lie bracket $[A,B]$, e.g.
	\[
	[T,R]=: {a\su{TR}W+b\su{TR}X+c\su{TR}T+d\su{TR}R}.
	\] 
\end{notation}

\begin{definition}
 Let $\sD=\ker\alpha\wedge\beta$ an Engel structure, we say that a framing $\sD=\bra W,\,X\ket$ is \emph{\abadapted} if $\sW=\bra W\ket$ spans the characteristic foliation, $c\su{WX}=\beta([W,X])=1$, and $d\su{XT}=\alpha([X,T])=1$.
\end{definition}
For given $\sD=\ker\alpha\wedge\beta=\bra W,\, X\ket$ it is always possible to rescale $W$ and $X$ to get a \abadapted\  framing.

The definition of $T$ and $R$ provides symmetries in the coefficients of the Lie brackets of the vector fields of the framing $\{W,\, X,\, T,\, R\}$.
The following result summarizes the ones that we will use in the rest of the paper.

\begin{lemma}\label{LEM_JacobiMadness}
	Let $\sD=\ker\alpha\wedge\beta$ be Engel and fix a \abadapted\  framing $\sD=\bra W,\, X\ket$.
	We have $c\su{WT}=c\su{WR}=c\su{XT}=c\su{XR}=0$ and $d\su{WX}=d\su{WT}=0$.
	Moreover
	\begin{subequations}
		\begin{align}
		b\su{WX}&=d\su{WR}\label{J1a}\\
		b\su{XT}&=-a\su{WT}\label{J1b}\\
		d\su{TR}&=\sL_Wd\su{XR}-\sL_Xd\su{WR}-a\su{WX}d\su{WR}-d\su{XR}d\su{WR}\label{J2a}\\
		c\su{TR}&=a\su{WR}+b\su{XR}\label{J2b}\\
		b\su{WR}&=-\sL_Wd\su{TR}+\sL_Td\su{WR}+a\su{WT}d\su{WR}+b\su{WT}d\su{XR}\label{J3a}\\
		b\su{TR}&=-\sL_Wc\su{TR}+d\su{WR}c\su{TR}\label{J3b}\\
		c\su{TR}&=-\sL_Xd\su{TR}+\sL_Td\su{XR}-a\su{WT}d\su{XR}+a\su{XT}d\su{WR}-b\su{XR}\label{J4a}\\
		a\su{TR}&=\sL_Xc\su{TR}-d\su{XR}c\su{TR} \label{J4b}
		\end{align}
	\end{subequations}

	\begin{proof}
		Equations $c\su{WT}=c\su{WR}=c\su{XT}=c\su{XR}=0$ are direct consequences of Definition~\ref{EQN_ReebDistribution}.
		Moreover $d\su{WX}=d\su{WT}=0$ follow from $\sL_W\sE\subset\sE$.
		Equations~\eqref{J1a}-\eqref{J4b} instead follow from all possible instances of the Jacobi identity.
		We prove $\eqref{J2a}$ and $\eqref{J2b}$, the other formulas follow similarly 
		\begin{align*}
		0&=\big[W,[X,R]\big]+\big[R,[W,X]\big]+\big[X,[R,W]\big]\\
		&=[W,a\su{XR}W+b\su{XR}X+d\su{XR}R]+[R,a\su{WX}W+b\su{WX}X+T]\\
		&+[X,-a\su{WR}W-b\su{WR}X-d\su{WR}R].
		\end{align*}
		We continue the calculation mod $\sD$
		\begin{align*}		 
		0&=b\su{XR}[W,X]+(\sL_Wd\su{XR})R+d\su{XR}[W,R]\\
		&-a\su{WX}[W,R]-b\su{WX}[X,R]-[T,R]\\
		&+a\su{WR}[W,X]-(\sL_Xd\su{WR})R-d\su{WR}[X,R]\qquad\mod\sD\\
		&=\Big(a\su{WR}+b\su{XR}-c\su{TR}\Big)T\\\
		&+\Big(\sL_Wd\su{XR}-\sL_Xd\su{WR}-a\su{WX}d\su{WR}-d\su{XR}d\su{WR}-d\su{TR}\Big)R.
		\end{align*}
	\end{proof}
	
\end{lemma}
Notice that Equations~\eqref{J3b} and~\eqref{J4b} provide another proof of Proposition~\ref{PROP_IntegrabilityOfR}.

\begin{remark}
	For a fixed defining form $\sE=\ker\alpha$ and any $X\in\Gamma\sD$ transverse to $\sW$ the form $\beta:=-\sL_X\alpha=-i_Xd\alpha$ is a defining form for $\sD$.
	Indeed, since $i_Wd\alpha=0$ on $\ker\alpha$, we have $\beta(W)=0$ and, by maximal non-integrability there is a section $Y\in\Gamma\sE$ such that $\beta(Y)\ne0$.
	This choice ensures that $d\su{XR}=\alpha([X,R])=-i_Xd\alpha(R)=\beta(R)=0$.
	In fact by a change $\beta\mapsto\beta+\nu\alpha$ we have complete freedom on the choice of $d\su{XR}$.
	Moreover up to rescaling $X$ we can ensure that $\bra W,\,X\ket$ is \abadapted.
\end{remark}

\section{Existence of integrable $\sR$}\label{SEC_IntegrableR}
Understanding whether a given Engel structures admits Engel defining forms that induce an integrable Reeb distribution is a complicated problem.
A first step in this direction is to compute the integrability condition~\eqref{EQN_integrabilityCondR} when we change defining forms.

\begin{lemma}\label{LEM_symmetriesOfForms}
	Let $\sD=\ker\alpha\wedge\beta$ be Engel and fix a \abadapted\ framing $\sD=\bra W,\, X\ket$.
	For another choice of Engel defining forms $\tilde\alpha$ and $\tilde\beta$, denote by $\tilde\sR=\bra\tilde T,\tilde R\ket$ the Reeb distribution associated with the new forms and $\tilde c\su{TR}=d\tilde\beta(\tilde R,\tilde T)$.
	We have
	\begin{enumerate}[1.]
		\item if $\tilde\alpha=\lambda\alpha$ and $\tilde\beta=\beta$ for $\lambda\in\sC^\infty(M)$ nowhere vanishing then
		\[
		\tilde R=\frac{1}{\lambda}R,\qquad \tilde T=T,\qquad \tilde c\su{TR}=\frac{1}{\lambda}c\su{TR}; 
		\]
		\item if $\tilde\alpha=\alpha$ and $\tilde\beta=\mu\beta$ for $\mu\in\sC^\infty(M)$ nowhere vanishing then
		\begin{align*}
		\tilde R&=R,\quad \tilde T=\frac{1}{\mu}\left(-\Big(\sL_X(\ln\mu)\Big)W+\Big(\sL_W(\ln\mu)\Big)X+T \right),\\
		\tilde c\su{TR}&=c\su{TR}+\sL_R(\ln\mu);
		\end{align*}
		\item if $\tilde\alpha=\alpha$ and $\tilde\beta=\nu\alpha+\beta$ for $\nu\in\sC^\infty(M)$ then
		\begin{align*}
		\tilde R&=\Big(-\nu^2-\sL_X\nu+\nu d\su{XR}\Big)W+\Big(\sL_W\nu-\nu d\su{WR}\Big)X-\nu T+R,\\
		\tilde T&=\nu W+T,\\
		\tilde c\su{TR}&=c\su{TR}-\nu\sL_W\nu-\sL_T(\nu)+\nu^2d\su{WR}+\nu d\su{TR}   
		\end{align*}
	\end{enumerate}
	
	\begin{proof}
		The proof of point 1 is straightforward.
		Let us prove point 2, the proof of 3 follows from a similar calculation; moreover these are the only points that we need in what follows.
		Since $\alpha$ is not changed and $\tilde\beta\wedge d\tilde\beta=\mu^2\beta\wedge d\beta$ we conclude $\tilde R=R$.
		Similarly we must have $\mu\tilde T=aW+bX+T$.
		Imposing $i_{\tilde T}(\tilde\alpha\wedge d\tilde\beta)=0$ yields the formula for $\tilde T$, here we need to use the hypothesis.
		The last formula follows directly from the evaluation $\tilde c\su{TR}=d\tilde\beta(\tilde T,\tilde R)$.
	\end{proof}
\end{lemma}

The previous formulas do not give any geometrical insight in understanding integrability of $\sR$.
We will instead turn to some stronger conditions that imply Equation~\eqref{EQN_integrabilityCondR}.
Notice that $d\beta^2=0$ happens if and only $c\su{TR}=d\beta(R,T)=0$, in particular this implies integrabiliy of $\sR$.
Proposition~\ref{PROP_PropertiesOfVolumePreservingW} ensuras that this is equivalent to the fact that $T$ preserves $\ker\beta$.
\begin{lemma}\label{LEM_dbeta2ConditionForIntegrability}
	For a given Engel structure $\sD=\ker\alpha\wedge\beta$ there exists a multiple $\tilde\beta=\mu\beta$ such that $d\tilde\beta^2=0$ if and only if for some $W$ and $X$ we have $a_{_{WR}}+b_{_{XR}}=0$.

	\begin{proof}
		Suppose that $a\su{WR}+b\su{XR}=0$.
		Notice that we cannot apply Lemma~\ref{LEM_JacobiMadness} because $\bra W,\,X\ket$ is not in general \abadapted.
		Nonetheless a similar calculation yields
		\begin{align*}
		0&=\beta\Big(\big[W,[X,R]\big]+\big[R,[W,X]\big]+\big[X,[R,W]\big]\Big)\\
		&=\beta\Big([W,a\su{XR}W+b\su{XR}X+d\su{XR}R]\Big)+\beta\Big([R,a\su{WX}W+b\su{WX}X+c\su{WX}T]\Big)\\
		&+\beta\Big([X,-a\su{WR}W-b\su{WR}X-d\su{WR}R]\Big)\\
		&=\beta\Big(b\su{XR}[W,X]+a\su{WR}[W,X]+(\sL_Rc\su{WX})T+c\su{WX}[R,T]\Big)\\
		&=\sL_Rc\su{WX}+c\su{WX}d\beta(T,R),
		\end{align*}
		where in the last equality we used the hypothesis.
		By maximal non-integrability we must have that $c\su{WX}=\beta([W,X])$  is nowhere vanishing, hence we can choose $\mu=c\su{WX}^{-1}$.
		Using the point 2 in Lemma~\ref{LEM_symmetriesOfForms} and the fact that $i_Rd\beta$ vanishes on $\sD$ we get
		\begin{align*}
		d(\mu\beta)(\tilde T,\tilde R)&=\Big(d\mu\wedge\beta+\mu d\beta\Big)\left(\frac{1}{\mu}T,R\right)\\
		&=-\frac{1}{\mu}\sL_R\mu+d\beta(T,R)=\frac{1}{c\su{WX}}\Big(\sL_Rc\su{WX}+c\su{WX}d\beta(T,R)\Big)=0.
		\end{align*}
		
		Conversely suppose that $d\tilde\beta^2=0$ and choose $\bra W,\,X\ket$ to be \abadapted.
		Lemma~\ref{LEM_symmetriesOfForms} ensures $\tilde R=R$ so that $\tilde a\su{WR}=a\su{WR}$ and $\tilde b\su{XR}=b\su{XR}$.
		We conclude by Equation~\eqref{J2b}.
	\end{proof} 
\end{lemma}

\subsection{Contactization}
We will now point out how to construct a contact structure starting from an Engel structure.
This construction is well-known to experts, but the author could not find any explicit reference.
\begin{definition}
	Let $(M,\, \sD=\ker\alpha\wedge\beta)$ be an Engel structure.
	The \emph{contactization of} $M$ is the contact $5$-manifold $(X=M\times \R,\, \xi=\ker\eta)$ with $\eta=\beta+s\alpha$, where we use $s$ for the coordinate along $\R$.
\end{definition}

The previous definition depends on the choice of $\alpha$ and $\beta$.
On the other hand we are interested in properties of $\xi$ which are invariant up to rescaling and translating the $\R$ factor.
The following result ensures that $\xi$ does not essentially dependent on the choice of Engel defining forms.
\begin{lemma}\label{LEM_VerificationsForContactization}
	Let $(M,\, \sD=\ker\alpha\wedge\beta)$ be an Engel structure and $X=M\times\R$ with $s$ indicating the coordinate on the $\R$-factor.
	The form $\eta=\beta+s\alpha$ defines a contact structure $\xi=\ker\eta$ on $X$.
	
	Moreover if we change Engel defining forms $\alpha\mapsto\tilde\alpha$ and $\beta\mapsto\tilde\beta$, there is a contactomorphism $\psi:(X,\,\ker\eta)\to (X,\,\ker\tilde\eta)$  of the form $\psi(p,\,s)=(p,\,f(p)s+g(p))$ where $f,\,g\in\sC^\infty(M)$ with $f$ nowhere vanishing.
	
	\begin{proof}
		To verify that $\eta$ is a contact form we calculate $d\eta=d\beta+sd\alpha+ds\wedge\alpha$ and
		\begin{equation}\label{EQN_deta2contactization}
		d\eta^2=d\beta^2+s^2d\alpha^2+2sd\alpha\wedge d\beta+2ds\wedge\alpha\wedge d\beta+2sds\wedge\alpha\wedge d\alpha.
		\end{equation}
		This, together with $\alpha\wedge\beta\wedge d\beta\ne0$ and $\alpha\wedge d\alpha\wedge\beta=0$, ensures that $\eta\wedge d\eta^2\ne0$.
		
		All possible choices of Engel defining forms can be written as $\tilde\alpha=\lambda\alpha$ and $\tilde\beta=\nu\alpha+\mu\beta$ for $\lambda,\,\nu,\,\mu\in\sC^\infty(M)$ with $\lambda$ and $\mu$ nowhere vanishing.
		Hence we have $\tilde\eta=\mu\beta+\nu\alpha+s\lambda\alpha$, and if we define 
		\[
		\psi(p,s):=\left(p,\,\frac{1}{\lambda(p)}\big(\mu(p)s-\nu(p)\big)\right)
		\]
		we have $\psi^*\tilde\eta=\mu\eta$.

	\end{proof}
\end{lemma}

Equation~\eqref{EQN_deta2contactization} implies that the Reeb vector field associated with $\eta$ takes the form (we use the notation introduced in~\ref{NOT_Commutators})
\begin{equation}\label{EQN_ReebForContactization}
R_\eta=T+sW+(c\su{TR}+sd\su{TR}+s^2d\su{WR})\partial_s.
\end{equation}
Since $d\beta^2=0$ if and only if $c\su{TR}=0$, Equation~\eqref{EQN_ReebForContactization} means that this happens if and only if $M\times\{0\}$ is invariant with respect to the Reeb flow.
\begin{remark}
	Consider the change of Engel defining forms $\tilde\alpha=\alpha$ and $\tilde\beta=\nu\alpha+\beta$.
	Lemma~\ref{LEM_VerificationsForContactization} ensures that there is a hypersurface $\tilde M\imm X$ graphical on $M\times\{0\}$ such that $(X,\,\tilde\eta)$ is (strictly contactomorphic to) the contactization of $(\tilde M,\,\sD=\ker\tilde\alpha\wedge\tilde\beta)$.
	Namely $\tilde M$ is the graph of $\nu:M\to\R$.
	
	Conversely take a hypersurface $\tilde M\imm X$ which is the graph of $g\in\sC^\infty(M)$, and define $\psi(p,s)=(p,\,s+g(p))$.
	The forms $\tilde\beta=\restrto{\eta}{\tilde M}=\beta+g\alpha$ and $\tilde\alpha=\restrto{(\sL_{\partial_s}\eta)}{\tilde M}=\alpha$ are defining forms for the Engel structure on $\tilde M$ obtained by pushing forward $\ker\alpha\wedge\beta$ via $\restrto\psi M:M\to\tilde M$.
\end{remark}

The previous remark, Remark~\ref{REM_EngelasIntersOfEvenContact} and Equation~\eqref{EQN_ReebForContactization} immediately give the proof of the following result.
\begin{proposition}
	Let $(M,\sD=\sE\cap\sE')$ be an Engel manifold, with $\sE$ and $\sE'$ even contact structures such that $\sE=[\sD,\sD]$, $\sW\subset\sE$ and $\sW'\pitchfork\sE$.
	There is a 1-to-1 correspondence between choices of $\sE'$ as above and hypersurfaces of $M\times\R$ which are graphical over $M$.
	
	Moreover, there is a choice of Engel defining forms $\sD=\ker\alpha\wedge\beta$ such that $d\beta^2=0$ if and only if there is a graphical hypersurface on the contactization $(X,\xi)$ invariant with respect to the Reeb flow associated with $\eta=\mu(\beta+s\alpha)$, for some $\mu\in\sC^\infty(M)$ nowhere vanishing.
\end{proposition}
The previous result gives a geometric interpretation of $d\beta^2=0$.
Unfortunately it is not very useful for practical purposes, because the dynamics of the Reeb vector field associated with $\eta$ can be very complicated.

\section{When are $\sD$ and $\sR$ totally geodesic?}\label{SEC_MetricProperties}
We now turn to the study of Riemannian properties of Engel structures.
Since a choice of Engel defining forms $\alpha$ and $\beta$ determines uniquely the splitting $TM=\sD\oplus\sR$, it is natural to consider Riemannian metrics $g$ which satisfy $\sD^\bot=\sR$.
In this context let $A$ and $B$ denote the $g$-duals of $\alpha$ and $\beta$ respectively, i.e.
\begin{equation}\label{EQN_DefinitionofABguys}
i_Ag=\alpha\spacedt{and}i_Bg=\beta.
\end{equation}

Since $\alpha$ and $\beta$ are linearly independent, the same is true for $A$ and $B$.
Moreover since $\sD^\bot=\sR$ they must be tangent to $\sR$, in particular $\sR=\bra A,B\ket$.
Equation~\eqref{EQN_DefinitionofABguys} implies
\[
\alpha(A)=\norm{A}^2,\quad \beta(B)=\norm{B}^2,\quad \alpha(B)=g(A,B)=\beta(A),
\]
and using these formulas we get
\begin{equation}\label{EQN_AandBasCombiofTandR}
A=g(A,B)T+\norm{A}^2R\spacedt{and} B=\norm{B}^2T+g(A,B)R.
\end{equation}

\begin{remark}
	In this context we have some freedom in the choice of $\alpha$.
	Suppose that $\alpha,\,\beta$, and $g$ are as above, and consider the new defining forms $\tilde\alpha=\lambda\alpha$ and $\tilde\beta=\beta$ for $\lambda\in\sC^\infty(M)$ nowhere vanishing.
	Lemma~\ref{LEM_symmetriesOfForms} ensures that the Reeb distribution $\tilde\sR$ associated with the new defining forms coincides with $\sR$.
	In particular, $g$ also satisfies $\sD^\bot=\tilde\sR$ and we have the formulas 
	\begin{equation}\label{EQN_SymmetriesOfAandB}
	\tilde A=\lambda A\spacedt{and}\tilde B=B.
	\end{equation}
\end{remark}

Recall that a distribution $\sH$ on a Riemannian manifold $(M,g)$ is \emph{totally geodesic} if for every $p\in M$ and $v_p\in\sH_p$ the geodesic through $p$ tangent to $v_p$ is tangent to $\sH$ at every point. 
We are interested in understanding under which conditions the distributions $\sD$ and $\sR$ are totally geodesic with respect to a metric $g$ making the splittin $TM=\sD\oplus\sR$ orthogonal.

\subsection{$\sR$ is never totally geodesic}
The goal of this section is to show that there is no metric $g$ such that $\sD^\bot=\sR$ and $\sR$ is totally geodesic.

\begin{lemma}\label{LEM_TotallyGeoCalculation}
	Let $\sD=\ker\alpha\wedge\beta$ be Engel and suppose that $g$ is a metric such that $\sD^\bot=\sR$, then $\sR$ is totally geodesic if and only if for all $U,\,U'\in\Gamma\sR$ and $V\in\Gamma\sD$ we have
	\[
	\sL_V(g(U,U'))+g([U',V],U)+g([U,V],U')=0.
	\]
	
	\begin{proof}
		It is well-known that $\sR$ is totally geodesic if and only if the following tensor vanishes (for the basic theory see \cite{rovenskii})
		\[
		h^\sR:\Gamma\sR\times\Gamma\sR\to\Gamma\sD\spacedt{s.t.}h^\sR(U,U')=\frac{1}{2}\sD\Big(\nabla_UU'+\nabla_{U'}U\Big),
		\]
		where $\sD$ and $\sR$ also denote, by abuse of notation, the orthogonal projections on the respective distributions.
		
		Now $h^\sR$ is zero if and only if for any $V\in\Gamma\sD$ we have $g(h^\sR(U,U'),V)=0$.
		Using Koszul's identity and the fact that $g(U,V)=0=g(U',V)$ we have
		\begin{align*}
		0=&4g(h^\sR(U,U'),V)=2g(\nabla_UU',V)+2g(\nabla_{U'}U,V)\\
		=&-\sL_V(g(U,U'))+g([U,U'],V)-g([U',V],U)+g([V,U],U')+\\
		&-\sL_V(g(U',U))+g([U',U],V)-g([U,V],U')+g([V,U'],U)\\
		=&-2\Big(\sL_V(g(U,U'))+g([U',V],U)+g([U,V],U')\Big).
		\end{align*} 
	\end{proof}
\end{lemma}

The following result furnishes an obstruction on the metric properties of the Reeb distribution associated with any pair of Engel defining forms, when the metric $g$ makes the splitting $TM=\sD\oplus\sR$ orthogonal.

\begin{proposition}
	Let $\sD=\ker\alpha\wedge\beta$ be an Engel structure and let $g$ be a metric such that $\sD^\bot=\sR$, then $\sR$ is not totally geodesic
	
	\begin{proof}
		Suppose that $\sR$ is totally geodesic and fix a framing $\sD=\bra W,\,X\ket$ which is \abadapted.
		Using Lemma~\ref{LEM_TotallyGeoCalculation} we get
		\begin{align*}
		g(h^\sR(A,A),W)=0&\ \RA&\sL_W\norm{A}^2+2\alpha([A,W])=0\\
		g(h^\sR(A,A),\,X)=0&\ \RA&\sL_X\norm{A}^2+2\alpha([A,X])=0\\
		g(h^\sR(B,B),W)=0&\ \RA&\sL_W\norm{B}^2+2\beta([B,W])=0\\
		g(h^\sR(B,B),\,X)=0&\ \RA&\sL_X\norm{B}^2+2\beta([B,X])=0\\
		g(h^\sR(A,B),W)=0&\ \RA&\sL_W(g(A,B))+\alpha([B,W])+\beta([A,W])=0\\
		g(h^\sR(A,B),\,X)=0&\ \RA&\sL_X(g(A,B))+\alpha([B,X])+\beta([A,X])=0\\
		\end{align*}
		Equation~\eqref{EQN_AandBasCombiofTandR} implies
		\begin{align*}
		[W,A]&=\Big(\sL_W(g(A,B))\Big)T+\Big(\sL_W\norm{A}^2+\norm{A}^2d\su{WR}\Big)R\\
		[W,B]&=\Big(\sL_W\norm{B}^2\Big)T+\Big(\sL_W(g(A,B))+g(A,B)d\su{WR}\Big)R\\
		[X,A]&=\Big(\sL_X(g(A,B))\Big)T+\Big(\sL_X\norm{A}^2+g(A,B)+\norm{A}^2d\su{XR}\Big)R\\
		[X,B]&=\Big(\sL_X\norm{B}^2\Big)T+\Big(\sL_X(g(A,B))+\norm{B}^2+g(A,B)d\su{XR}\Big)R.
		\end{align*}
		These in turns yield
		\begin{enumerate}[a.]
			\item $\sL_W\norm{A}^2+2\norm{A}^2d\su{WR}=0$
			\item $\sL_X\norm{A}^2+2g(A,B)+2\norm{A}^2d\su{XR}=0$
			\item $\sL_W\norm{B}^2=0$
			\item $\sL_X\norm{B}^2=0$
			\item $\sL_W(g(A,B))+g(A,B)d\su{WR}=0$
			\item $\sL_X(g(A,B))+\norm{B}^2+g(A,B)d\su{XR}=0$
		\end{enumerate}
		Now (e) and (d) imply that $\norm{B}^2$ is constant on orbits of $W$ and $X$.
		Since $\sD=\bra W,X\ket$ is bracket-generating, Chow's Theorem implies that $\norm{B}^2$ is constant.
		Notice that the same formulas must hold for $\tilde A$ and $\tilde B$ given by Equation~\eqref{EQN_SymmetriesOfAandB}. 
		Since $\|\tilde A\|=|\lambda|\norm{A}$, we can suppose that $\norm{A}^2$ is also constant, say
		\[
		\norm {A}=1.
		\]
		In particular for some $c\in\R$, we get
		
		\begin{enumerate}[a.]
			\item $d\su{WR}=0$
			\item $g(A,B)+d\su{XR}=0$
			\item[c./d.] $\norm{B}=c$
			\item[e.] $\sL_W(g(A,B))=0$
			\item[f.] $\sL_X(g(A,B))+\norm{B}^2+g(A,B)d\su{XR}=0$
		\end{enumerate}
		Hence (b) and (f) together give
		\[
		\sL_X(g(A,B))=g(A,B)^2-\norm{B}^2.
		\]
		The Cauchy-Schwarz inequality reads $g(A,B)^2\le\norm{A}^2\norm{B}^2$ with equality if and only if $\{A,\,B\}$ is linearly dependent.
		Since this never happens, the inequality is sharp and
		\[
		\sL_X(g(A,B))<\norm{A}^2\norm{B}^2-\norm{B}^2=0.
		\]
		This is a contradiction because we must have $\sL_X(g(A,B))=0$ on the critical points of the function $p\mapsto \restrto{g(A,B)}{p}$.
	\end{proof}
\end{proposition}

The geometric reason for this obstruction is not clear.
Notice that the hypothesis on $g$ cannot be relaxed as Example~\ref{EXM_TorusKEngel} furnishes an Engel structure $\sD=\ker\alpha\wedge\beta$ on $T^4$ where $\sR$ is a totally geodesic foliation with respect to the standard metric.
In this case the splitting $TM=\sD\oplus\sR$ is not orthogonal.

\subsection{$\sD$ totally geodesic implies $\sR$ integrable}
In this section we study the properties of the Reeb distribution associated with a totally geodesic Engel structure $\sD=\ker\alpha\wedge\beta$.
Let $g$ be such that $\sD^\bot=\sR$ and choose a framing $\sD=\bra W,\,X\ket$ which is orthonormal and such that $W$ spans the characteristic foliation.
The proof of the following lemma is exactly the same as the proof of Lemma~\ref{LEM_TotallyGeoCalculation}.
\begin{lemma}\label{LEM_TotallyGeoEngelIFF}
	Under the above hypothesis $\sD$ is totally geodesic if and only if for all $X_1,X_2\in\{W,X\}$ and $Y\in\Gamma\sR$ we have
	\[
	g([X_1,Y],X_2)+g([X_2,Y],X_1)=0
	\]
\end{lemma}
The previous result permits us to express the Lie brackets of sections of $\sD$ with sections of $\sR$ in a simple way.
\begin{corollary}\label{COR_TotallyGeoFraming}
	Under the above hypothesis $\sD$ is totally geodesic if and only if
	\begin{align*}
	[W,T]&=b\su{WT}X&[W,R]&=b\su{WR}X+d\su{WR}R\\
	[X,T]&=-b\su{WT}W+d\su{XT}R&[X,R]&=-b\su{WR}W+d\su{XR}R\\
	\end{align*}
	
	\begin{proof}
		Since $i_Td\beta=-c\su{TR}\alpha$, all $T$-components must vanish.
		Moreover $[W,T]\in\sE$ hence $d\su{WT}=0$.
		Hence the only components left to calculate are the ones in the direction of $W$ and $X$.
		Since $\{W,\,X\}$ is an orthonormal basis and $\sD^\bot=\sR$ we can calculate them using Lemma~\ref{LEM_TotallyGeoEngelIFF}:
		\[
		a\su{WT}=g([W,T],W)=0,\quad a\su{WR}=g([W,R],W)=0,
		\]
		similarly $b\su{XT}=0$ and $b\su{XR}=0$.
		Moreover
		\[
		b\su{WT}=g([W,T],X)=-g([X,T],W)=-a\su{XT}
		\]
		and $b\su{WR}=-a\su{XR}$, which concludes the proof.
	\end{proof}
\end{corollary}

The following result links metric properties of $\sD$ with integrability properties of $\sR$.
\begin{corollary}
	Let $\sD=\ker\alpha\wedge\beta$ be Engel and let $g$ be a metric such that $\sD^\bot=\sR$.
	If $\sD$ is totally geodesic, then there exists a nowhere vanishing function $\mu\in\sC^\infty(M)$ such that $\tilde\beta=\mu\beta$ satisfies $d\tilde\beta^2=0$.
	In particular, $\tilde\sR$ is integrable.
	\begin{proof}
		Corollary~\ref{COR_TotallyGeoFraming} ensures that $a\su{WR}+b\su{XR}=0$, which is exactly the hypothesis of Lemma~\ref{LEM_dbeta2ConditionForIntegrability}.
	\end{proof}
\end{corollary}

The converse of this result is likely false, but we do not have a counterexample.
The structure on $Nil^4$ constructed in Section~\ref{SEC_GeoEngel} is an example of totally geodesic Engel structure with respect to a compatible metric.

\section{Engel vector fields}\label{SEC_EngelVF}

We now turn our attention to the study of Engel structures that admit symmetries.
The existence of $1$-parameter families of contactomorphism for any given contact structure is well-known.
For Engel structures on the other hand the existence of such families of symmetries is tightly related to the dynamics of the characteristic foliation.

\begin{definition}
	Let $(M,\sD)$ be an Engel structure an \emph{Engel vector field} is a vector field whose flow preserves $\sD$.
\end{definition}

\begin{remark}\label{REM_EngelVFPreservesTheFlag}
	If $Z$ preserves $\sD$ then automatically it must preserve its Engel flag, i.e.
	\[
	\sL_Z\sW=\sW,\quad\sL_Z\sD=\sD\spacedt{and}\sL_Z\sE=\sE.
	\]
	In \cite{montgomeryPaper} there is an example of Engel structure admitting a unique $1$-parameter family of symmetries.
	It is unclear if there are Engel structure which do not admit any 1-parameter families of symmetries.
\end{remark}

The following result furnishes a relation between existence of symmetries of $\sD$ and the dynamics of $\sW$.
\begin{lemma}\label{LEM_TransverseSymmetryRAVolPres}
    Let $\sD$ be an Engel structure and suppose that $Z$ is an Engel vector field transverse to $\sE$.
    Then there exists a pair of defining forms $\sD=\ker\alpha\wedge\beta$ such that $Z=R$ and $d\alpha^2=0$.
	\begin{proof}
		Use Proposition~\ref{PROP_PropertiesOfVolumePreservingW} to find $\alpha$ such that $\alpha(Z)=1$ and $d\alpha^2=0$.
		We need to find $\beta$ so that $Z=R$.
		
		By Remark~\ref{REM_EngelVFPreservesTheFlag} the flow of $Z$ must preserve the Engel flag $\sW\subset\sD\subset\sE$.
		This means that we can choose a framing $\sE=\bra W,\,X,\,Y\ket$ satisfying 
		\begin{equation}\label{EQN_ImmaOnlyBeUsedOnce}
		\sL_Z\bra W\ket=\bra W\ket,\quad\sL_Z\bra X\ket\subset\bra W,X\ket,\quad\sL_Z\bra Y\ket\subset\bra W,X,Y\ket.
		\end{equation}
		Choose $\beta$ so that $\ker\beta=\bra W,\,X,\,Z\ket$ and $\beta(Y)=1$.
		Exactly as in the proof of Lemma~\ref{LEM_ExistenceOfDefiningForms} the forms $\alpha$ and $\beta$ are Engel defining forms for $\sD$.
		Now Equation~\eqref{EQN_ImmaOnlyBeUsedOnce} implies $[W,Z]=aW$, so that 
		\[
		d\beta(Z,W)=\beta([W,Z])=a\beta(W)=0.
		\]
		Similarly we have $d\beta(Z,X)=0$, so that $d\beta=0$ on $\ker\beta$.
		This implies $Z\in\ker\beta\wedge d\beta$, and since $\alpha(Z)=1$ we conclude $Z=R$. 
	\end{proof}
\end{lemma}

Proposition~\ref{PROP_PropertiesOfVolumePreservingW} furnishes various equivalent conditions to the existence of a transverse even contact symmetry.
The following result gives an adapted pair of Engel defining forms if such a symmetry exists for $\sE=[\sD,\sD]$.

\begin{lemma}\label{LEM_EquivCondToL_RD=D}
	Let $(M,\sD)$ be an Engel structure trivial as a bundle and $M$ orientable.
	The following are equivalent
	\begin{enumerate}[1.]
		\item $\sW$ is volume preserving;
		\item there exist $\alpha$ and $\beta$ such that $\sL_R\sD\subset\sD$;
		\item there exist $\alpha$ and $\beta$ such that $\ker d\alpha=\bra W,\,R\ket$;
		\item there exist $\alpha$ and $\beta$ such that $\beta\wedge d\alpha=0$;
		\item there exists $\alpha$ such that the conformal class of $\beta=\sL_X\alpha$ does not depend on the choice of $X\in\Gamma\sD$ transverse to $\sW$.
	\end{enumerate}
	Moreover if $d\alpha^2=0$ there is a choice of $\beta=-\sL_X\alpha$ such that all the above properties are verified simultaneously.
	\begin{proof}
		\fbox{$1\RA2$} This is a corollary of Proposition~\ref{PROP_PropertiesOfVolumePreservingW} and Lemma~\ref{LEM_TransverseSymmetryRAVolPres}.
		
		\fbox{$2\RA3$} The hypothesis implies that $\sL_R\sE\subset\sE$ so that $\sL_R\alpha=i_Rd\alpha=\lambda\alpha$ for some $\lambda\in\sC^\infty(M)$.
		Since $i_W(\alpha\wedge d\alpha)=0$ we have $i_Wd\alpha=h_W\alpha$ for some $h_W\in\sC^\infty(M)$.
		Hence $d\alpha^2=0$ if and only if $(i_Wd\alpha)(R)=0$, but we have $(i_Wd\alpha)(R)=-(i_Rd\alpha)(W)=-\lambda\alpha(W)=0$.
		This also proves that $\ker d\alpha=\bra W,R\ket$.
		
		\fbox{$3\RA4$} This is obvious since both $W$ and $R$ are in the kernel of $\beta$ and of $d\alpha$.
		
		\fbox{$4\RA5$} For any choice of $X$ we must have $0=i_X(\beta\wedge d\alpha)=-\beta\wedge \sL_X\alpha$, which is only possible if $\sL_X\alpha$ is a multiple of $\beta$.
		
		\fbox{$5\RA1$} By Proposition~\ref{PROP_PropertiesOfVolumePreservingW} it suffices to proof that $i_Wd\alpha=0$.
		Suppose this is not true, then as in the above proof we must have $i_Wd\alpha=h_W\alpha$ with $h_W\in\sC^\infty(M)$ not identically zero.
		For any given $X\in\Gamma\sD$ transverse to $W$, the vector field $\tilde X=X+W$ is transverse to $W$, but $\beta=i_Xd\alpha$ is not a multiple of $\tilde\beta=i_{\tilde X}d\alpha=\beta+h_W\alpha$.  
		
		The last statement follows directly from this proof.
	\end{proof}
\end{lemma}

\section{Engel Killing vector fields}\label{SEC_EngelKillingVF}
In this section we will study Engel structures admitting transverse Killing symmetries.
\begin{proposition}\label{PROP_ExistenceOfEKimpliesGoodForms}
	Let $(M,\sD)$ be an Engel structure and let $g$ be a Riemannian metric.
	Suppose that $Z\in\fX(M)$ is Engel, Killing and orthogonal to $\sE$, then there exists a choice of defining forms $\alpha$ and $\beta$ such that $d\alpha^2=0=d\beta^2$ and $\beta\wedge d\alpha$.
	\begin{proof}
		Since $Z\bot\ker\alpha$ it must be in particular transverse to it, hence Lemma~\ref{LEM_TransverseSymmetryRAVolPres} implies the existence of $\alpha$ and $\beta$ such that $R=Z$ and $d\alpha^2=0$.
		Up to rescaling, we can suppose that $\norm{R}=1$.
		Fix an orthonormal basis $\sD=\bra W,X\ket$ and complete it with a vector field $Y$ to an orthonormal basis of $\sE$.
		This implies that $\{W,\,X,\,Y,\,R\}$ is an orthonormal framing.
		
		Since $R$ is an Engel vector field, as in the proof of Lemma~\ref{LEM_TransverseSymmetryRAVolPres}, we have
		\begin{align*}
		[W,R]&=a\su{WR}W\\
		[X,R]&=a\su{XR}W+b\su{XR}X\\
		[Y,R]&=a\su{YR}W+b\su{YR}X+c\su{YR}Y.
		\end{align*}
		Since $R$ is Killing we have
		\[
		0=(\sL_Rg)(W,W)=\sL_R(g(W,W))-2g(\sL_RW,W)=2a\su{WR},
		\]
		so that $[W,R]=0$.
		Similarly $b\su{XR}=c\su{YR}=0$.
		Moreover
		\[
		0=(\sL_Rg)(W,X)=\sL_R(g(W,X))-g(\sL_RW,X)-g(W,\sL_RX)=a\su{XR},
		\]
		so that $[X,R]=0$.
		Similarly $a\su{YR}=0$ and $b\su{YR}=0$.
		Hence we have $[W,R]=[X,R]=[Y,R]=0$.
		If we now pick $\beta=-\sL_X\alpha$ we immediately have $d\beta^2=0$ and $\beta\wedge d\alpha=0$, so that $\beta=-\sL_X\alpha$ follows from Lemma~\ref{LEM_EquivCondToL_RD=D}.
	\end{proof}
\end{proposition}

\begin{definition}
	A \emph{K-Engel structure} is a triple $(M,\,g,\,Z)$ where $\sD$ is an Engel structure, $g$ is a metric and $Z$ is a vector field which is Engel, Killing and orthogonal to $\sE$.
\end{definition}
Moreover the Engel defining forms $\alpha$ and $\beta$ satisfying $d\alpha^2=0=d\beta^2$ and $d\alpha\wedge\beta=0$ are called \emph{K-Engel forms}.

\begin{corollary}\label{COR_KEngelFraming}
	Let $(M,\,\sD,\,g,\,Z)$ be a K-Engel structure, then there exists Engel defining forms $\alpha$ and $\beta$ and a \abadapted\  framing $\sD=\bra W,\,X\ket$ such that $Z=R$ and we have
	\begin{align*}
	[W,R]&=[X,R]=[T,R]=0\\
	[W,X]&=a\su{WX}W+T\\
	[W,T]&=a\su{WT}W+b\su{WT}X\\
	[X,T]&=a\su{XT}W-a\su{WT}X+R
	\end{align*}
	where the functions $a\su{WX},\,a\su{WT},\,b\su{WT}$ and $a\su{XT}$ are constant on the orbits of $R$.
	\begin{proof}
		The proof of Proposition~\ref{PROP_ExistenceOfEKimpliesGoodForms} implies the existence of a framing $\{W,\,X,\,Y,\,R\}$ such that $R$ commutes with every vector field in the framing.
		Now $\alpha\wedge d\beta=0$ implies $d\alpha\wedge d\beta=0$, which in turn implies $[T,R]=0$.
		
		We need to rescale $W$ and $X$ so that $c\su{WX}=1=d\su{XT}$. 
		This is possible because the Jacobi identity implies
		\begin{align*}
		\sL_R\Big(\beta([W,X])\Big)&=\big(\sL_R\beta\big)\Big([W,X]\Big)-\beta\Big(\big[R,[W,X]\big]\Big)\\
		&=\beta\Big(\big[X,[R,W]\big]+\big[W,[X,R]\big]\Big)=0,
		\end{align*}
		and similarly $\sL_R(\alpha([X,T]))=0$.
		Hence we can rescale $W$ and $X$ as follows
		\[
		W\mapsto\frac{\alpha([X,Y])}{\beta([W,X])}W,\qquad X\mapsto\frac{1}{\alpha([X,Y])}X
		\]
		to get a new framing of $\sD$ which satisfies all previous conditions and is \abadapted.
		We have
		\begin{align*}
		[W,R]&=[X,R]=[T,R]=0\\
		[W,X]&=a\su{WX}W+b\su{WX}X+T\\
		[W,T]&=a\su{WT}W+b\su{WT}X\\
		[X,T]&=a\su{XT}W+b\su{XT}X+R.
		\end{align*}
		Equation~\eqref{J1a} implies $b\su{WX}=d\su{WR}=0$ and Equation~\eqref{J1b} implies $b\su{XT}=-a\su{WT}$.
	\end{proof}
\end{corollary}

We will often denote only by $\sD=\ker\alpha\wedge\beta$ the K-Engel structure $(\sD=\ker\alpha\wedge\beta,\,g,\,Z=R)$, if we do not want to put an accent on the metric $g$.
We call \emph{K-Engel framing} the framing $TM=\bra W,\,X,\,T,\,R\ket$ defined in the previous corollary.

\begin{remark}
 The choice of K-Engel defining forms and framing is not unique.
 Let $(M,\,g,\,Z)$ is a K-Engel structure and fix K-Engel forms $\alpha$ and $\beta$.
 All other possible choices of K-Engel defining forms are $\tilde\alpha=\alpha$ and $\tilde\beta=\mu\beta$ with $\mu\in\sC^\infty(M)$ constant on $R$-orbits and nowhere-vanishing.
 
 Moreover if $\{W,\,X,\,T,\,R\}$ is a K-Engel framing all other choices of K-Engel framings are $\tilde W=\lambda W$, $\tilde X=\mu X+\nu W$, $\tilde T$ given by Lemma~\ref{LEM_symmetriesOfForms}, and $\tilde R=R$, where $\lambda,\mu,\nu\in\sC^\infty(M)$ are constant on $R$-orbits and with $\lambda$ and $\mu$ nowhere-vanishing.

%
 \end{remark}

The converse of Proposition~\ref{PROP_ExistenceOfEKimpliesGoodForms} is not true.
The existence of defining forms satisfying $d\alpha^2=0=d\beta^2$ and $d\alpha\wedge\beta=0$ only ensures that $R$ is a Killing vector field if $\sL_R$ acts in a diagonalizable way on $\sD$.
\begin{proposition}\label{PROP_CondOnFormsRAKEngel}
	Let $\sD=\ker\alpha\wedge\beta$ be an Engel structure such that $d\alpha^2=0=d\beta^2$ and $d\alpha\wedge\beta=0$.
	Suppose that there exists $X\in\Gamma\sD$ transverse to $\sW$ and such that $\sL_R\bra X\ket=\bra X\ket$ for $b\in\sC^\infty(M)$.
	Then there exists a a metric $g$ such that $(\sD,\,g,\,Z)$ is K-Engel.
	
	\begin{proof}
		The idea is to construct a framing $\{W,\,X,\,T,\,R\}$ such that $R$ commutes with every vector field in it, and then take the metric $g$ making this framing orthonormal.

		First of all notice that Lemma~\ref{LEM_EquivCondToL_RD=D} ensures that the flow of $R$ preserves $\sD$ and its Engel flag.
		In particular for any framing $\sD=\bra W,\,X\ket$ we must have that $\sL_RW=aW$ and $\sL_RX=bX+cW$ for some smooth functions $a,\,b,$ and $c$.
		The hypothesis ensures that we can choose $X$ such that $c=0$.
		Up to rescaling we can suppose that $\bra W,\,X\ket$ is \abadapted.
		
		The condition $d\beta^2=0$ implies that $[T,R]$ is a multiple of $R$.
		Since $d\alpha\wedge\beta=0$ implies $d\alpha\wedge d\beta=0$, we must have $[T,R]=0$.
		Using $\sL_R\beta=0$ we get
		\begin{align*}
		0&=\sL_R\Big(\beta\big([W,X]\big)\Big)=-\beta\Big(\big[R,[W,X]\big]\Big)\\
		&=\beta\Big(\big[X,[R,W]\big]+\big[W,[X,R]\big]\Big)=-a-b.
		\end{align*}
		Similarly $\sL_R\alpha=0$ and $[T,R]=0$ imply $0=\sL_R\Big(\alpha\big([X,T]\big)\Big)=-b$.
		Hence the claim.
	\end{proof}
\end{proposition}

\section{Some remarks on the dual of $W$}\label{SEC_rho}
Before continuing with the discussion on general K-Engel structures we make some observations about some special cases.
If $(M,\sD)$ is Engel and its flag are orientable, a choice of a framing of $\sD=\bra W,\,X\ket$, with $\sW=\bra W\ket$ and of defining forms $\alpha$ and $\beta$ yields the framing $\{W,\,X,\,T,\,R\}$.
We are interested in the properties of the dual coframing $\{\rho,\,\tau,\,\beta,\,\alpha\}$ in the case where $d\alpha^2=0$ and $\beta=-\sL_X\alpha$.
Under these hypothesis we can determine whether $\alpha$ and $\beta$ are K-Engel by looking at $\rho$.
\begin{lemma}
	Suppose $\sD=\ker\alpha\wedge\beta$ satisfies $d\alpha^2=0$ and $\beta=-\sL_X\alpha$.
	Then $\alpha$ and $\beta$ are K-Engel if and only if there is a choice of $W$ and $X$ such that $\sL_R\rho=0$ modulo $\bra\beta\ket$.
	\begin{proof}
		The key observation is that $\sL_R\rho=-a\su{WR}\rho-a\su{XR}\tau$ modulo $\bra\beta\ket$.
		If $\alpha$ and $\beta$ are K-Engel then Corollary~\ref{COR_KEngelFraming} ensures that there is a choice of $W$ and $X$ such that they commute with $R$.
		This implies in particular $a\su{WR}=0$ and $a\su{XR}=0$.
		
		Conversely suppose we have such framing.
		Up to rescaling $X$, we have $c\su{WX}=1$ and this does not change $\rho$ and $R$.
		Notice that
		\[
		d\su{XT}=d\alpha(T,X)=\beta(T)=1,
		\]
		so that $\bra W,\, X\ket$ is \abadapted.
		The hypothesis $d\alpha^2=0$ and $\beta=-\sL_X\alpha$ imply $d\su{WR}=d\su{XR}=d\su{TR}=0$.
		These together with Equation~\eqref{J4a} imply $c\su{TR}=-b\su{XR}$.
		Equation~\eqref{J2b} reads $c\su{TR}=a\su{WR}+b\su{XR}$, so using $\sL_R\rho\wedge\beta=0$ yields
		\[
		2c\su{TR}=a\su{WR}=0,
		\]
		which translates to $d\beta^2=0$.
		Now all hypothesis of Proposition~\ref{PROP_CondOnFormsRAKEngel} except possibly for $\sL_RX=bX$, but this is again a consequence of $\sL_R\rho\wedge\beta=0$, since this means $a\su{XR}=0$.
		
	\end{proof}
	
\end{lemma}

Notice that the hypothesis of the previous lemma is verified if $d\rho=0$.
Some aspects of the theory of K-Engel structures such that $d\rho=0$ resemble the theory of Sasakian manifolds (see \cite{blair} Section 6.8).
The proof of the following is borrowed by the analogous result in the Sasakian case.
I wish to thank Giovanni Placini for explaining it to me.
\begin{theorem}\label{THM_CupLengthKEngel}
	Let $\sD=\ker\alpha\wedge\beta$ be K-Engel, fix the induced K-Engel framing $\{W,\,X,\,T,\,R\}$ and its dual coframing $\{\rho,\,\tau,\,\beta,\,\alpha\}$.
	If $d\rho=0$ then the cup-length of $M$ is smaller than $4$.
	\begin{proof}
		The idea is to prove that $R$ is contained in the kernel of any harmonic $1$-form.
		In this case for any $a_1,...,a_4\in H^1(M,\R)$ we pick harmonic representatives $\theta_1,...,\theta_4$ and the cup product will be the class of $\theta_1\wedge\cdots\wedge\theta_4$, but this is zero because $R$ is in its kernel.
		
		Since $R$ is a Killing vector field, its flow acts by isometries, hence it sends harmonic forms to harmonic forms.
		Moreover it acts trivially in cohomology.
		These two facts imply $\sL_R\theta=0$ for all $\theta$ harmonic 1-form.
		Write $\theta=\eta+f\alpha$ where $\eta(R)=0$ and $f=\theta(R)$, we want to prove that $f=0$.
		Using the fact that $\theta$ is closed we have
		\[
		0=\sL_R\theta=d(\theta(R))=df,
		\]
		hence $f$ is constant and $0=d\theta=d\eta+fd\alpha$.
		A simple calculation yields
		\[
		d(\rho\wedge\eta\wedge\alpha)=d\rho\wedge\eta\wedge\alpha-\rho\wedge d\eta\wedge\alpha+\rho\wedge\eta\wedge d\alpha=
		-f\rho\wedge\alpha\wedge d\alpha,
		\]
		where we used $d\rho=0$ and the fact that $\rho\wedge\eta\wedge d\alpha=0$ because $R$ is in the kernel of all the forms (see Lemma~\ref{LEM_EquivCondToL_RD=D}).
		Now $\rho\wedge\alpha\wedge d\alpha$ is a volume form because $\rho(W)=1$, hence if $f\ne0$ this would imply that it is also exact, which is impossible since $M$ is closed.

	\end{proof}
	
\end{theorem}

Example~\ref{EXM_TorusKEngel} furnishes K-Engel structures on $T^4$, so the hypothesis $d\rho=0$ in the previous lemma is crucial.

\section{Topology of K-Engel manifolds}\label{SEC_TopologyOfKEngel}
We will now focus on topological obstructions to the existence of a K-Engel structure.
Corollary~\ref{COR_KEngelFraming} ensures that there exists a framing $\{ W,\,X,\,T,\,R\}$ such that $R$ commutes with all vector fields in the framing.
This fact has strong consequences on the topology of $M$, which come from the theory of transverse structures to a foliation.
Results on transversally parallelizable foliations in \cite{molino} permit to prove that $M$ has a fibre bundle structure, where the fibres are the closure of the orbits of $R$.
In our context we can prove a stronger result, namely that $M$ has the structure of a principal torus bundle.
The torus acting on $M$ will be the closure of the flow of $R$ in the isometry group of $(M,g)$.
The first step is the following result, which is a consequence of Chow's Theorem.
\begin{lemma}\label{LEM_OrbitsOfRKEngel}
	Let $(M,\,\sD=\ker\alpha\wedge\beta)$ be a K-Engel and $p,q\in M$, then there exists an isotopy $\psi_t:M\to M$ which commutes with the flow of $R$ and such that $\psi_1(p)=q$.
	In particular all $R$ orbits are equivariantly isotopic on $M$. 
	\begin{proof}
		By hypothesis $R$ is transverse to $\ker\alpha$ and we have a framing $\{W,\,X,\,T\}$ of vector fields commuting with $R$.
		This means that if we can join $p$ and $q$ with a piecewise smooth path $\gamma$ obtained by glueing together pieces of orbits of $W,\,X$ and $T$ we will obtain the isotopy $\psi_t$ by integrating the path.
		The existence of such a path is exactly the statement of Chow's Theorem (see \cite{gromov}).
	\end{proof}
\end{lemma}

We can use Lemma~\ref{LEM_OrbitsOfRKEngel} to describe the topology of $M$.
\begin{lemma}\label{LEM_PrincipalBundleStructure}
 Let $(M,\,\sD=\ker\alpha\wedge\beta)$ be a K-Engel structure, then the closure of the orbits of $R$ are all tori $T^k$.
 Moreover $M$ is a principal $T^k$-bundle over a manifold of dimension $4-k$.
 In this case we say that $(M,\,\sD=\ker\alpha\wedge\beta)$ has rank $k$.
 \begin{proof}
   Let $\phi_t^R$ be the flow of $R$, this is a $1$-parameter subgroup of the compact Lie group $G$ of isometries of $(M,g)$.
   Since it is Abelian, its closure is an embedded torus $T^k\subset G$.
   This means that the action of $T^k$ on $M$ is effective so that the intersection of all its stabilizers is the trivial group.
   Since $T^k$ is Abelian, the principal orbit type must be the trivial one (for a proof see Example I.2.6 in \cite{audin}).
   This ensures that there exists a point $p\in M$ whose orbit is an embedded $T^k$.
   Now by Lemma~\ref{LEM_OrbitsOfRKEngel} all orbits have the same type, so they are all embedded copies of $T^k$.
   In particular the action is free and $M$ is a principal $T^k$-bundle.
 \end{proof}
\end{lemma}

We now discuss all possible cases for $k\in\{1,2,3,4\}$.
If $k=1$ then $R$ is totally periodic and this case plays a central role, as ensured by the following result.
\begin{proposition}\label{PROP_KEngelDeformToRegular}
 Let $(M,\,\sD=\ker\alpha\wedge\beta,\,R)$ be a K-Engel structure of rank $k$, then there exist K-Engel structures of rank $1$ $(M,\,\sD=\ker\alpha_i\wedge\beta_i,\, R_i)$ for $i=1,...,k$ such that $R_1,...,R_k$ are linearly independent.
 \begin{proof}
  As above let $T^k$ be the closure of $\phi^R_\R$ in the group of isometries of $(M,g)$.
  Fix a basis $\bra A_1,...,A_k\ket$ for the Lie algebra of $T^k$ and consider the vector fields $R_i=\exp(tA_i)$ for $i=1,...,k$.
  We can pick $A_i$ so that $\alpha(R_i)$ is nowhere-vanishing for all $i=1,...,k$.
  Moreover by construction these are totally periodic vector fields whose flows preserve $\alpha$ and $\beta$ and a K-Engel framing $\{W,\,X,\,T,\,R\}$.
  Now we claim that
  \[
   \alpha_i=\frac{1}{\alpha(R_i)}\alpha\spacedt{and}\beta_i=-\sL_X\alpha_i
  \]
  are K-Engel defining forms for $i=1,...,k$.
  First of all using $\sL_R\alpha=0$, $\sL_R(\alpha(R_i))=0$ and $\alpha_i(R_i)=1$ we get
  \[
   i_{R_i}d\alpha_i=\sL_{R_i}\alpha_i=0,
  \]
  so that $d\alpha_i^2=0$.
  By definition $\beta_i\wedge d\alpha_i=0$, moreover using $[R_i,X]=0$ we get $\sL_{R_i}\beta_i=0$, which in turn implies $d\beta_i^2=0$.
  Finally Proposition~\ref{PROP_CondOnFormsRAKEngel} ensures that $(M,\,\sD=\ker\alpha_i\wedge\beta_i,\, R_i)$ is K-Engel.
 \end{proof}

\end{proposition}

The previous result is the analogue of Theorem 7.1.10 in \cite{boyer} which asserts that every K-contact structure can be perturbed to a quasi-regular one. 
Since they play such a central role, we will analyse the properties of K-Engel structures of rank $1$ in Section~\ref{SEC_EngelBW}.
In Section~\ref{SEC_T2Bundles} we will give some further constraints on the bundle structure for K-Engel of rank $2$.

Suppose now that $k=3$, this implies that $M$ is a $T^3$-principal bundle over $S^1$.
Since the classifying space for principal $T^3$-bundles is simply connected, this must be the trivial bundle, i.e. $M=T^4$.
The last case is $k=4$, we prove by contradiction that this cannot happen.
In this situation we have that $M=T^4$ and the orbits of $R$ are dense.
The conditions $\sL_R\alpha=\sL_R\beta=0$ imply that these forms are homogeneous, which in turn means that they must be closed, which contradicts non-integrability.

The previous discussion furnishes a proof of the following result.
\begin{theorem}\label{THM_TopologyOfKEngel}
	If $(M,\sD)$ admits a K-Engel structure, then $M$ is diffeomorphic to one of the following:
	\begin{itemize}
		\item $T^4$;
		\item a principal $T^2$-bundle over a surface;
		\item a principal $S^1$-bundle over a $3$-manifold.
	\end{itemize} 
\end{theorem}

It is unclear which principal torus bundles admit K-Engel structures.
In Sections~\ref{SEC_T2Bundles} and~\ref{SEC_EngelBW} we will characterise them via some differential conditions, which are nonetheless hard to verify directly.
We end this section with the construction of a family of K-Engel structures on $T^4$ providing examples for which the rank takes all possible values.
\begin{example}\label{EXM_TorusKEngel}
	On $\R^4$ with coordinates $(x,y,z,t)$ consider the distribution given by 
	\[
	\sD=\Big\bra W=\cos(2\pi t)\partial_x+\sin(2\pi t)\partial_y+\partial_z,\, X=\partial_t\Big\ket.
	\]
	This defines an Engel structure (in fact this is the Lorentz prolongation of the standard Lorentz structure on $\R_{2,1}$).
	We choose defining forms
	\[
	\alpha=dz-\cos(2\pi t)dx-\sin(2\pi t)dy,
	\qquad
	\beta=-\sin(2\pi t)dx+\cos(2\pi t)dy.
	\]	
	An explicit calculation yields
	\[
	R=\partial_z
	\spacedt{and}
	T= -\sin(2\pi t)\partial_x+\cos(2\pi t)\partial_y.
	\]
	Since $R$ is Engel and Killing for the metric making $\{W,\,X,\,T,\,R\}$ orthonormal, so that this is a K-Engel structure.
	
	Up to choosing a lattice $\Lambda$ in $\R^4$ in the right way we can make sure that this structure passes to the quotient $T^4=\R^4/\Lambda$.
	Moreover we can control the dimension of the closure of the orbits of $R$.
	The only condition on $\Lambda=\bra e_1,...,e_4\ket$ is that the $t$-components of the vectors must be integers (otherwise $\cos(2\pi t)$ and $\sin(2\pi t)$ will not pass to the quotient).
	Picking $e_4=(0,0,0,1)$ leaves complete freedom for $e_1,\, e_2$ and $e_3$ in the orthogonal space to $e_4$.
	Since $R=\partial_z$ we can make sure that the closure of its orbits in the quotient is $S^1,\, T^2$, or $T^3$.
	This gives a family of K-Engel structures $\sD_s$ on $T^4$ smoothly depending on $s$ and such that the rank assumes all possible values as $s$ changes.
\end{example}

\subsection{K-Engel $T^2$-bundles}\label{SEC_T2Bundles}
In this section we will study K-Engel structures of rank $2$.
This implies that $M$ is a principal $T^2$-bundle over a surface $\Sigma_g$ of genus $g$.
If not otherwise specified we will always denote by $R_1$ and $R_2$ a choice of two totally periodic vector fields whose flows give a splitting of $S^1\times S^1=T^2$ acting on $M$.
Moreover if defining forms $\sD=\ker\alpha\wedge\beta$ are fixed, we can suppose, as in the proof of Proposition~\ref{PROP_KEngelDeformToRegular}, that $R_1$ and $R_2$ are both transverse to $\sE$ and $R=\epsilon R_1+R_2$ for some irrational $0<\epsilon<1$.

It is well known that principal $T^2$-bundles are classified up to isomorphism by their \emph{Euler class} (see \cite{walczak}).
Let $\pi:M\to\Sigma_g$ be the bundle projection, and let $\Omega$ be a generator of $H^2(\Sigma_g,\Z)$.
The Euler class of the the $T^2$-bundle $e(\pi)\in H^2(\Sigma_g,\Z)\oplus H^2(\Sigma_g,\Z)$ can be identified with a pair of integers $(n_1,n_2)$.
Quotienting $M$ via the action of $R_2$ (resp. $R_1$) we get a (oriented) $3$-manifold $N_1$ (resp $N_2$), which, in turn, is an $S^1$-bundles over $\Sigma_g$ with Euler class $n_1\Omega$ (resp. $n_2\Omega$).
Hence we get the commuting diagram
\[
\begin{tikzcd}[column sep=small]
& M\arrow[dl,"q_2"']\arrow[dr,"q_1"]\arrow[dd,"\pi"] & \\
N_1\arrow[dr,"p_1"'] & & N_2\arrow[dl,"p_2"]\\
&\ \Sigma_g&
\end{tikzcd}
\]
In particular $p_i:N_i\to\Sigma_g$ has Euler class $n_i\Omega$ and $M\to N_i$ has Euler class $n_2p_i^*\Omega$.

It is unclear which  $T^2$-bundle admits K-Engel structures of rank $2$.
The following result gives a (rather obscure) differential condition which characterises such bundles.

\begin{lemma}
 Let $\Sigma_g$ be a surface, $\Omega$ and orientation form on it, and $n_1,n_2\in\Z$.
 Consider a principal $T^2$-bundle over $(\Sigma_g,\Omega)$ with Euler class $(n_1,n_2)$ and fix connection forms $(\theta_1,\theta_2)$ such that $\theta_i(R_j)=\delta_{ij}$.
 Then $M$ admits a K-Engel structure $\sD=\ker\alpha\wedge\beta$ of rank $2$ if and only if there exist two functions $f,g\in\sC^\infty(\Sigma_g)$ and two $1$-forms $\alpha_0,\beta_0\in\Omega^1(\Sigma_g)$ satisfying
 \begin{equation}\label{EQN_ConditionFor2KEngel}
  \begin{cases}
    df(p)\ne0\spacedt{or}(f(p)N+n_2)\Omega(p)+d\alpha_0(p)\ne0\quad\forall p\in\Sigma_g\\
   Ng^2\,\Omega+g\,d\beta_0+\beta_0\wedge dg\ne0\\
   Ng\,\Omega+g\,d\alpha_0+\beta_0\wedge df=0
  \end{cases}
 \end{equation}
 where $N=n_1-\epsilon n_2$ and $R=\epsilon R_1+R_2$.
 
 In this case we have Engel defining forms
 \begin{equation}\label{EQN_2KEngelDefiningForm}
    \begin{split}
    \alpha&=\pi^*f\ \theta_1+(1-\epsilon\, \pi^*f)\theta_2+\pi^*\alpha_0\\
    \beta&=\pi^*g\ (\theta_1+\epsilon\theta_2)+\pi^*\beta_0.
    \end{split}
 \end{equation}

\begin{proof}
  Suppose that $M$ admits a K-Engel structure $(M,\,\sD=\ker\alpha\wedge\beta)$ of rank $2$.
  By hypothesis we have $R=\epsilon R_1+R_2$ and we set $\tilde f=\alpha(R_1)$ and $\tilde g=\beta(R_1)$.
  These functions are invariant so there must exist $f,g\in\sC^\infty(\Sigma_g)$ such that $\tilde f=\pi^*f$ and $\tilde g=\pi^*g$.
  Since $\alpha$ and $\beta$ are $R$-invariant, for some $1$-forms $\alpha_0,\beta_0\in\Omega^1(\Sigma_g)$ we have
  \begin{align*}
    \alpha&=\pi^*f\ \theta_1+(1-\epsilon\, \pi^*f)\theta_2+\pi^*\alpha_0\\
    \beta&=\pi^*g\ (\theta_1+\epsilon\theta_2)+\pi^*\beta_0.
  \end{align*}
  The Engel conditions
  \[
   \begin{cases}
    \alpha\wedge d\alpha\ne0\\
    \alpha\wedge\beta\wedge d\beta\ne0\\
    \alpha\wedge d\alpha\wedge\beta=0
   \end{cases}
  \]
  translate to (using $d\theta_i=n_i\pi^*\Omega$)
  \[
   \begin{cases}
    \pi^*df\wedge\theta_1\wedge\theta_2+\pi^*\Big(f\tilde\Omega+\alpha_0\wedge df\Big)\wedge\theta_1
    +\pi^*\Big((1-\epsilon\,f)\tilde\Omega+\alpha_0\wedge df\Big)\wedge\theta_2\ne0\\
    \pi^*(Ng^2\,\Omega+g\,d\beta_0+\beta_0\wedge dg)\wedge\theta_1\wedge\theta_2\ne0\\
    \pi^*(Ng\,\Omega+g\,d\alpha_0+\beta_0\wedge df)\wedge\theta_1\wedge\theta_2=0
   \end{cases}
  \]
  where $\tilde\Omega=(fN+n_2)\Omega+d\alpha_0$.
  Notice that the first equation is verified if and only if for every point $p\in\Sigma_g$ on which $df(p)=0$ we have $\tilde\Omega(p)\ne0$, since $f(p)\ne0$.
  Hence we get Equation~\eqref{EQN_ConditionFor2KEngel}.
  
  Viceversa suppose that we have $f,g$ and $\alpha_0,\beta_0$ such that Equation~\eqref{EQN_ConditionFor2KEngel} is verified.
  The same calculation as above ensures that the defining forms $\alpha$ and $\beta$ given by Equation~\eqref{EQN_2KEngelDefiningForm} are Engel defining forms.
  Moreover by construction they satisfy $d\alpha^2=0=d\beta^2$ and $\beta\wedge d\alpha=0$ and $R=\epsilon R_1+R_2$, so that it preserves any invariant framing $\bra W,\, X\ket$ of $\sD$.
  Proposition~\ref{PROP_CondOnFormsRAKEngel} ensures that $\sD=\ker\alpha\wedge\beta$ is a K-Engel structure.
\end{proof}
\end{lemma}

\subsection{Engel Boothby-Wang}\label{SEC_EngelBW}
In this section we study properties of K-Engel structures of rank $1$.
The construction that we will present has already appeared, in a different context, in the work of Mitsumatsu \cite{mitsu} under the name of \emph{prequantum prolongation}.

Let $\sD=\ker\alpha\wedge\beta$ be K-Engel of rank $1$, so that $M$ is an $S^1$-bundle $\pi:M\to N$ over a $3$-manifold $N$ with $R$ tangent to fibres.
Since $\sL_R\alpha=0$ and $\alpha(R)=1$ we have that $\alpha$ is a connection form.
This means that $d\alpha$ descends to a closed $2$-form $\omega$ satisfying $[\omega]\in H^2(N,\Z)$ (because it represents the Euler class of the bundle). 

Notice that $\sL_R\beta=0$ so that, since we have fixed a connection, we get a $1$-form $\lambda$ on $N$ by pushing down $\beta$.
Since $\ker\beta\wedge d\beta=\bra R\ket$, the form $\lambda$ is a contact form.
Similarly $W$ descends to a Legendrian vector field $L=\pi_*W$ which is divergence free (with respect to the contact volume) and such that $i_L\omega=0$.
The goal of this section is to reverse this construction.

\begin{remark}
	Let $(N,\lambda)$ be a contact structure on a $3$-manifold, and let $L$ be a non-singular Legendrian vector field.
	Then $L$ is the kernel of a closed non-singular $2$-form $\omega$ if and only if a rescaling of $L$ preserves the contact volume $\lambda\wedge d\lambda$.
\end{remark}

Let $N$ be a closed $3$-manifold and $p:H^2(N,\Z)\to H^2(N,\R)$.
Recall that for any closed integral $2$-form $\omega$ there is an $S^1$-bundle $\pi:M\to N$ whose Euler class maps to $[\omega]$ via $p$.
If there is torsion in $H^2(N,\Z)$ there are finitely many choices for the isomorphism type of $M$, any of these would work.

The space of connections on $\pi:M\to N$ is an affine space on the space of closed $1$-forms on $N$.
We fix an arbitrary choice of connection $\alpha_0$.
If two closed forms $\theta_1$ and $\theta_2$ differ by an exact form, then there is a gauge transformation isotopic to the identity that sends the connection form $\alpha_0+\theta_1$ to $\alpha_0+\theta_2$.
This means that the moduli space of possibile choices of connections is (at most) an affine space on $H^1(N,\R)$.

\begin{proposition}\label{PROP_EngelBW}
	Let $(N^3,\lambda)$ be a contact structure, $L$ be a Legendrian vector field such that its dual $2$-form $\omega=i_L(\lambda\wedge d\lambda)$ is integral and closed, and fix $[\theta]\in H^1(N,\R)$.
	Then the principal $S^1$-bundle $\pi:M\to N$ with Euler class $[\omega]$ admits a K-Engel structure $\sD=\ker\alpha\wedge\beta$, such that $\alpha$ is the connection form associated to $[\theta]$ and $\beta=\pi^*\lambda$.
	
	\begin{proof}
		First of all notice that $\alpha$ is uniquely determined by the discussion above.
		We need to verify that $\alpha$ and $\beta$ are K-Engel defining forms.
		Now $\alpha$ defines an even contact structure because by definition $d\alpha=\pi^*\omega$, so that $\alpha\wedge d\alpha\ne0$.
		
		Since $\lambda$ is contact $\beta$ is even contact.
		Moreover $\ker\beta\wedge d\beta$ is spanned by the vector field $R$ tangent to the fibres and normalized by $\alpha$.
		This implies in turn that $\alpha\wedge\beta\wedge d\beta\ne0$.
		Finally $\alpha\wedge d\alpha\wedge\beta=0$, since already $\beta\wedge d\alpha=\pi^*(\lambda\wedge\omega)=0$ because $i_L(\lambda\wedge\omega)=0$.
		So $\sD$ is an Engel structure.
		
		Now for dimensional reasons $d\alpha^2=0=d\beta^2$.
		As we have already seen $d\alpha\wedge\beta=0$ and $R$ acts on $\sD$ in a diagonalizable way, as can be seen by choosing $X=\pi^*\tilde L$ where $\tilde L$ is a Legendrian line field nowhere tangent to $L$.
		Proposition~\ref{PROP_CondOnFormsRAKEngel} implies that $\sD$ is K-Engel.
	\end{proof}
\end{proposition}

\begin{remark}
  We refer to the construction in Proposition~\ref{PROP_EngelBW} as \emph{Engel Boothby-Wang construction}.
  This is the Engel analogue of the Boothby-Wang construction \cite{bw}.
  The above discussion proves that all K-Engel structures such that $R$ is totally periodic are obtained via an Engel Boothby-Wang construction.
  Example~\ref{EXM_TorusKEngel} provides K-Engel structures which do not come from an Engel Boothby-Wang construction, because their rank is not $1$.
\end{remark}

We finish this section providing some examples of $4$-manifolds admitting K-Engel structures of rank $1$.
Section~\ref{SEC_GeoEngel} provides further examples.
\begin{example}
	We will now show that both Cartan and Lorentz prolongations (see \cite{cp3}) can admit K-Engel structures.
	Consider the surface $\Sigma_g$ equipped with a metric $h$ of constant scalar curvature $k$ and construct the unit circle bundle $S^1\to N\to\Sigma_g$.
	Denote by $A\in\fX(N)$ a unit vector field tangent to the fibres.
	The Levi-Civita connection induces a choice of horizontal bundle on $N$ and we denote by $B$ a tautological vector field, i.e. for $l\in T_pN$ of unit norm we want $\pi_*(B(p,l))\in \R l$.
	Finally we choose $C$ so that $\{B,C\}$ is an orthonormal basis for the horizontal bundle of $N$.
	It is a classical result that
	\[
	[A,B]=C,\qquad[B,C]=kA\spacedt{and}[C,A]=kB.
	\]
	Consider $M=N\times S^1$ and call $t$ the coordinate along $S^1$.
	
	One can pick K-Engel defining forms
	\[
	\alpha=\cos t\ a+\sin t\ b+c\spacedt{and}\beta=-\sin t\ a+\cos t\ b
	\]
	and verify that the following is a K-Engel framing
	\begin{align*}
	  W&=\cos t\ A+\sin t\ B-C+(k+1)\partial_t\\
	  X&=\partial_t\\
	  T&=-\sin t\ A+\cos t\ B\\
	  R&=C-k\partial_t.
	\end{align*}
	This is obtained as Lorentz prolongation of the the conformal structure on $N$ having $\{A,\,B,\,C\}$ as an orthonormal basis, with $A$ and $B$ negative and $C$ positive.
	
	Another K-Engel structure on this manifold is provided by the oriented Cartan prolongation of the contact structure $\bra A,\,B\ket$ on $N$.
	Namely take $\sD=\bra W,\, X\ket$ where $W=\partial_t$ and $X=\cos t\ A+\sin t\ B$.
	Then the even contact structure $\sE$ is spanned by $\sD$ and 
	\[
	  Y:=[W,X]=-\sin t\ A+\cos t\ B.
	\]
	If we take $R=C+k\partial_t$ we can verify that it commutes with $W,\,X$ and $Y$ so that $\{W,\,X,\,Y,\,R\}$ is a K-Engel framing and taking the dual coframing $\{\rho,\,\tau,\,\beta,\,\alpha\}$ we get K-Engel defining forms $\alpha$ and $\beta$.
\end{example}

\begin{example}
	If an orientable $T^2$-bundle over a surface admits an $S^1$-action tangent to the fibres then the monodromy $\rho:\pi_1(\Sigma_g)\to Diff^+(\Sigma_g)$ has the form
	\[
	\rho(\gamma)=\Matrix{1&\lambda(\gamma)\\0&1}
	\]
	for $\lambda(\gamma)\in\Z$ (see Proposition 4.4 in \cite{walczak}).
	
	Take the flat $T^2$-bundle over $T^2$ with monodromy given by $\rho(a_i)=A_i$ for $\pi_1(T^2)=\bra a_1,\,a_2\ket$, where
	\[
	A_i=\Matrix{1&\lambda_i\\0&1},\qquad\textup{with }\lambda_i\in\Z,\,i=1,2.
	\]
	Take coordinates $\tilde M=\R^2\times T^2=\{(x,y,u,v)\}$, then the $1$-forms
	\[
	a=du+(\lambda_1x+\lambda_2y)dv\spacedt{and}b=dv
	\]
	are invariant with respect to the transformations
	\[
	\left(\Matrix{x\\y},\,\Matrix{u\\v}\right)\mapsto\left(\Matrix{x-1\\y},\,A_1\Matrix{u\\v}\right)
	\]
	and
	\[
	\left(\Matrix{x\\y},\,\Matrix{u\\v}\right)\mapsto\left(\Matrix{x\\y-1},\,A_2\Matrix{u\\v}\right)
	\]
	hence they define $1$-forms on $M$.
	Similarly the formulas 
	\[
	U=\partial_u
	\spacedt{and}
	V=\partial_v-(\lambda_1x+\lambda_2y)\partial_u
	\]
	define nowhere-vanishing vector fields tangent to the fibres of $M$.
	The forms
	\[
	\alpha=a-\cos v\,dz-\sin v\,dy
	\spacedt{and}
	\beta=-\sin v\,dx+\sin v\,dy
	\]
	are K-Engel forms with $R=U$.
\end{example}

\section{Contact fillings}\label{SEC_ContactFillings}
There is a way to see Engel structures as special submanifolds of contact $5$-dimensional manifolds.
The K-Engel structures coming from the Engel Boothby-Wang construction are examples of such submanifolds in compact contact $5$-manifolds.

Let $(X^5,\,\xi=\ker\eta)$ be a contact structure on an orientable $5$-manifold $X$, and let $M$ be an orientable embedded hypersurface $M^4\to X^5$ transverse to $\xi$.
This means that (locally) we can find a Legendrian vector field $L\in\fX(X)$ transverse to $M$.
This data permits to define two $1$-forms on $M$ as follows
\begin{equation}\label{EQN_alphaandbetaForEngelTypeHypersurfaces}
\beta:=\restrto{\eta}{M}\spacedt{and}\alpha=\restrto{(\sL_L\eta)}{M}.
\end{equation}

We look for conditions such that a neighbourhood of $M$ is contactomorphic to the contactization of the Engel structure having defining forms $\alpha$ and $\beta$.
In what follows we will always suppose that $X$ is oriented by the contact volume and we will always take oriented embeddings $M\to X$ if not otherwise specified.
We have the following consequence of the Contact Weinstein's Neighbourhood Theorem.
\begin{lemma}\label{LEM_TransverseHypersurfaceNeighbourhood}
	Let $(X,\,\eta,\,L,\,M)$ be as above and define $\alpha$ and $\beta$ as in Equation~\eqref{EQN_alphaandbetaForEngelTypeHypersurfaces}.
	Then there is an open neighbourhood $\Op(M)=M\times(-\epsilon,\epsilon)$ and a function $f:M\times(-\epsilon,\epsilon)\to\R$ such that
	\[
	f\eta=\beta+s\alpha
	\]
	on $\Op(M)$, where we denote by $s$ the coordinate along $(-\epsilon,\epsilon)$.
	
	\begin{proof}
		Use the flow $\phi_t$ of $L$ to construct an embedding of a tubular neighbourhood $\psi:M\times(-\epsilon,\epsilon)\to X$ such that $\psi(p,s)=\phi_s(p)$.
		We identify $M\equiv M\times\{0\}$, $\alpha\equiv\psi^*\alpha$ and $\beta\equiv\psi^*\beta$, and we define $\eta_1=\beta+s\alpha$.
		We have $L=\partial_s$ and $\eta_0=\psi^*\eta$ is a contact form on $M\times(-\epsilon,\epsilon)$.
		We want to prove that $\eta_1$ is also a contact form on $\Op(M)$ and that
		\begin{equation}\label{EQN_guy1}
		\restrto{\eta_0}{M}=\restrto{\eta_1}{M}\spacedt{and}\restrto{d\eta_0}{M}=\restrto{d\eta_1}{M}
		\end{equation}
		so that the Contact Weinstein's Neighbourhood Theorem gives us a map $\tilde\psi:M\times(-\epsilon,\epsilon)\to M\times(-\epsilon,\epsilon)$ which satisfies $\tilde\psi^*\eta_1=f\eta$.
		
		Since $\phi_0=id$, Equation~\eqref{EQN_alphaandbetaForEngelTypeHypersurfaces} implies Equation~\eqref{EQN_guy1}.
		A direct calculation yields
		\[
		\eta_1\wedge d\eta_1^2=2\Big(ds\wedge\alpha\wedge\beta\wedge d\beta+sds\wedge\alpha\wedge\beta\wedge d\alpha\Big),
		\]
		hence if $\alpha\wedge\beta\wedge d\beta\ne0$ on $M$ we conclude that $\eta_1$ is contact on a (possibly smaller) tubular neighbourhood of $M$.
		Plugging $L$ into $\eta\wedge d\eta^2\ne0$ we obtain $\eta\wedge i_Ld\eta\wedge d\eta\ne0$.
		Since the kernel of $i_L(\eta\wedge d\eta^2)$ is $L$, and this is transverse to $M$, we have
		\[
		0\ne\restrto{\Big(i_L(\eta\wedge d\eta^2)\Big)}{M}=
		-2\restrto{\Big(\eta\wedge d\eta\wedge i_Ld\eta\Big)}{M}=
		-2\beta\wedge d\beta\wedge\alpha.
		\]
	\end{proof}
\end{lemma}

\begin{remark}\label{REM_IssueWithL}
	In the previous theorem we cannot ensure in general that $\tilde\psi_* L=\partial_s$.
	On the other hand we will only be interested in the quantities $\sL_L\eta$ and $\sL_Ld\eta$ on $M$, and the formula $f\eta=\beta+s\alpha$ implies that $\sL_{\partial_s}(f\eta)=\alpha=\sL_L\eta$ and $\sL_{\partial_s}d(f\eta)=d\alpha=\sL_Ld\eta$ on $M$.
\end{remark}

We say that $L$ \emph{preserves the contact volume on} $M$ if
\[
\restrto{\sL_L\left(\eta\wedge d\eta^2\right)}{p}=0\qquad\forall p\in M.
\]
This condition will be useful in proving that $\alpha$ and $\beta$ are Engel defining forms under certain hypothesis on $L$.
Moreover the following remark says that it is not too difficult to achieve.

\begin{remark}
	For any given $(X,\,\eta,\,L,\,M)$ as above, we can rescale $\eta$ so that $L$ preserves the contact volume on $M$.
	Let $g:X\to \R$ be such that $\sL_L\left(\eta\wedge d\eta^2\right)=g\,\textup{vol}_X$.
	Since $L$ is transverse to $M$, there is function $\lambda$ satisfying $3\restrto{\sL_L\lambda}{M}=-\restrto{g}{M}$.
	For every $p\in M$ we have
	\begin{align*}
	\restrto{\sL_L\left(e^\lambda\eta\wedge d(e^\lambda\eta)^2\right)}{p}&=
	\restrto{\sL_L\left(e^{3\lambda}\eta\wedge d\eta^2\right)}{p}\\
	&=\restrto{(3\sL_L\lambda)e^{3\lambda}\left(\eta\wedge d\eta^2\right)}{p}+e^{3\lambda} \restrto{\sL_L\left(\eta\wedge d\eta^2\right)}{p}=0.
	\end{align*}
\end{remark}

\begin{definition}
	Let $(X^5,\,\xi=\ker\eta)$ be a contact structure and $M\to N$ an embedded hypersurface.
	Let $L$ be a Legendrian line field transverse to $M$ which preserves the contact volume on $M$.
	We say that $M$ is an \emph{Engel-type hypersurface} if $\alpha=\restrto{(\sL_L\eta)}{M}$ is an even contact structure on $M$.
\end{definition}

The previous definition is justified by the following result.
\begin{lemma}
	Let $(X^5,\ker\eta)$ be a contact structure and $M$ be an Engel-type hypersurface, then $\alpha=\restrto{(\sL_L\eta)}{M}$ and $\beta=\restrto{\eta}{M}$ are Engel defining forms for $\sD=\ker\alpha\wedge\beta$.
	Moreover there exists a neighbourhood of $M$ in $X$ contactomorphic to the contactization of $(M,\sD)$.
	\begin{proof}
		By assumption we have $\alpha\wedge d\alpha\ne 0$.
		Lemma~\ref{LEM_TransverseHypersurfaceNeighbourhood} ensures that we have a neighbourhood of $M$ such that $f\eta=\beta+s\alpha$.
		Moreover the proof of the lemma ensures that $\alpha\wedge\beta\wedge d\beta\ne0$.
		
		The formula
		\[
		f\eta\wedge d(f\eta)^2=2\Big(ds\wedge\alpha\wedge\beta\wedge d\beta+sds\wedge\alpha\wedge\beta\wedge d\alpha\Big)
		\]
		ensures that
		\[
		\alpha\wedge\beta\wedge d\alpha=i_{\partial_s}\restrto{\sL_{\partial_s}\left(f\eta\wedge d(f\eta)^2\right)}{M}.
		\]
		Now since $\restrto{f\eta}{M}=\restrto{\eta}{M}$ and $\restrto{\sL_{\partial_s}(f\eta)}{M}=\restrto{\sL_L\eta}{M}$, by Remark~\ref{REM_IssueWithL} we have $\restrto{d(f\eta)}{M}=\restrto{d\eta}{M}$ and $\restrto{\sL_{\partial_s}(d(f\eta))}{M}=\restrto{\sL_Ld\eta}{M}$, hence
		\[
		\alpha\wedge\beta\wedge d\alpha=
		i_{\partial_s}\sL_{\partial_s}\restrto{\left(f\eta\wedge d(f\eta)\right)}{M}=
		i_{\partial_s}\sL_L\restrto{\left(\eta\wedge d\eta\right)}{M}=0.
		\]
	\end{proof}
\end{lemma}

\begin{example}
	All Engel manifolds $(M,\sD)$ appear as Engel-type hypersurfaces of their contactization.
	An example of Engel-type hypersurface on a compact contact manifold is provided by the Engel structure on $S^3\times S^1$ given by the Boothby-Wang construction on $(S^3,\,\lambda=a,\,\omega=db,\,[0])$.
	Here $a$ is the standard contact form on $S^3$ and $b$ is the $1$-form obtained from it by multipling by the quaternion $j$.
	In standard coordinates $\C^2=(z_1=x_1+iy_1,z_2=x_2+iy_2)$ we have
	\begin{align*}
	a&=-y_1dx_1+x_1dy_1-y_2dx_2+x_2dy_2\\
	b&=y_2dx_1+x_2dy_1-y_1dx_2-x_1dy_2.
	\end{align*}
	The Engel Boothby-Wang construction yields an Engel structure on $M=S^3\times S^1$ with defining forms $\alpha=dt+b$ and $\beta=c$.
	We have a contact structure on $X=N\times D^2$ which is given by the kernel of $\eta=\beta+r^2\alpha$ where $r$ is the radial coordinate.
	Then $M=\partial X$ is an Engel-type hypersurface with $rL=\partial_r$.
\end{example}

The previous example and the analogous definitions for the symplectic case motivates the following definition.
\begin{definition}
	A \emph{contact filling} of an Engel structure $(M,\sD)$ is a contact manifold $(X,\eta)$ such that $M=\partial X$ is an Engel-type hypersurface.
\end{definition}

It is unclear which Engel structures admit a contact filling.
The following result ensures that Engel Boothby-Wang manifolds are fillable.
\begin{theorem}
	Let $(M,\sD)$ be obtained from $(N,\,\lambda,\,\omega,\,[\theta])$ via the Boothby-Wang construction.
	Then $M$ admits a contact filling.
	\begin{proof}
		By hypothesis $\pi:M\to N$ is an $S^1$-bundle and the Engel defining forms are a connection form $\alpha$ and $\beta=\pi^*\lambda$.
		Consider now the disc bundle $X\to N$ of $M$ and denote by $r$ the radial coordinate on each fibre.
		The form $\eta=\beta+r^2\alpha$ is contact and $M=\partial X$ is an Engel-type hypersurface, with $rL=\partial_r$.
	\end{proof}
	
\end{theorem}

\section{Geometric K-Engel structures}\label{SEC_GeoEngel}
We end the paper providing the complete list of geometric K-Engel structures.
A geometric Engel structure is an Engel structure $\sD$ on a manifold $M$ admitting a Thurston geometry $(X,G)$ and such $\sD$ is invariant with respect to the action of $G$.
\begin{theorem}[\cite{vogelGeometricEngel}]\label{THM_GeoEngel}
 The following Thurston geometries admit a geometric Engel structure which is unique up to equivalence
 \begin{itemize}
  \item $(S^3\times\R,\,SU(2)\times\R)$
  \item $(\widetilde{Sl}(2,\R)\times\R,\,\widetilde{Sl}(2,\R)\times\R)$
  \item $(Nil^3\times\R,\,Nil^3\ltimes\R)$
  \item $(Sol^4(m,n),\,Sol^4(m,n))$
  \item $(Sol_0^4,\,Sol^4(\lambda))$
  \item $(Sol_1^4,\,Sol_1^4)$
  \item $(Nil^4,\,Nil^4)$
 \end{itemize}
No subgeometry of the other Thurston geometries of dimension $4$ admits geometric Engel structures.
\end{theorem}
For a proof and a more detailed description of the geometries see \cite{vogelGeometricEngel}.
We say that a K-Engel structure $(\sD,\,g,\,Z)$ is geometric if $\sD$, $g$, and $Z$ are left-invariant.
If one can construct a left-invariant framing $\{W,\,X,\,Y,\,R\}$ such that $\sW=\bra W\ket$, $\sD=\bra W,\,X\ket$, $\sE=\bra W,\,X,\,Y\ket$ and $R$ commutes with all other vector fields in the framing, then $G$ admits a geometric K-Engel structure.
Indeed, if this happens, we can take the metric $g$ which makes the framing orthonormal, and $(\sD,\,g,\,R)$ is K-Engel.
We can construct K-Engel forms $\alpha$ and $\beta$ by considering the dual coframing $\{\rho,\,\tau,\,\beta,\,\alpha\}$.

\begin{remark}\label{REM_GeometricKEngel}
 The converse of the above construction is also true.
 Indeed suppose that $(\sD,\,g,\,R)$ is a geometric K-Engel structure and consider a left-invariant orthonormal framing $\sD=\bra W,\,X\ket$, where $\sW=\bra W\ket$.
 A calculation similar to the one in the proof of Proposition~\ref{PROP_ExistenceOfEKimpliesGoodForms} implies that $R$ commutes with both $W$ and $X$.
 This implies that we can construct left-invariant K-Engel forms and framing.
 
 Notice in particular that if $\sD$ is left-invariant, and it admits a section $X$ that does not commute with any linearly independent left-invariant vector field $R$, then $\sD$ cannot admit a geometric K-Engel structure.
\end{remark}

We will now go through the list in Theorem~\ref{THM_GeoEngel} and verify which of these Engel structures admit geometric K-Engel structures.
In what follows we refer to Section 3 in \cite{vogelGeometricEngel} for the missing details.

\subsection{$S^3\times\R,\,\widetilde{Sl}(2,\R)$}
We can treat the first two geometries in Theorem~\ref{THM_GeoEngel} simultaneously since their Lie Algebra is generated by $A,\,B,\,C$, and $\partial_t$ such that the only non-vanishing Lie brackets are
\[
 [A,B]=C\qquad[B,C]=A\qquad[C,A]=kB.
\]
where $k=1$ for $S^3\times\R$, and $k=-1$ for $\widetilde{Sl}(2,\R)$.
Here $\partial_t$ is tangent to the $\R$-factor.

We have the Engel structure $\sD=\bra W,\, X\ket$ where $W=\partial_t+B$ and $X=A$.
Then
\[
 Y:=[W,X]=-C,\qquad Z:=[X,Y]=kB,
\]
so that $R=\partial_t$ commutes with $W,\,X$, and $Y$, and $\{W,\,X,\,Y,\,R\}$ satisfies the desired properties.

\subsection{$Nil^3\times\R$}
In the case of $Nil^3\times\R$ there is no geometric Engel structure invariant with respect to the full isometry group.
We instead have to consider the subgeometry with isometry group $G$ whose Lie algebra is generated by $A,\,B,\,C$, and $D$ satisfying
\[
 [A,B]=C\qquad[D,A]=-B\qquad[D,B]=A,
\]
and such that all other Lie brackets vanish.
In this case we can take $\sD=\bra W,\, X\ket$ where $W=D$ and $X=A$, so that
\[
 Y:=[W,X]=-B,\qquad Z:=[X,Y]=-C.
\]
Now $R=Z$ commutes with $W,\,X$, and $Y$, so that $\{W,\,X,\,Y,\,R\}$ satisfies the desired properties.

\subsection{$Sol^4(m,n)$}
In this case the Lie algebra is generated by $X_1,\,X_2,\,X_3$, and $T$, such that the only non-vanishing Lie brackets are $[T,X_i]=c_iX_i$ for $i=1,2,3$, where $c_i$ are real numbers satisfying $c_1+c_2+c_3=0$, $c_1>c_2>c_3$, and such that $e^{c_i}$ are the roots of
\begin{equation}\label{EQN_Pmn}
  P(m,n)=-\lambda^3+m\lambda^2-n\lambda+1.
\end{equation}
Up to isomorphism the only geometric Engel structure on $Sol^4(m,n)$ is $\sD=\bra W,\,X\ket$, where $W=X_1+X_2+X_3$ and $X=T$.

Suppose that a left-invariant vector field $Z=a_1X_1+a_2X_2+a_3X_3+aT$ commutes with $X$, this implies $a_ic_i=0$ for $i=1,2,3$, hence $X$ and $Z$ can only be linearly independent if $c_i=0$ for some $i$.
Hence if $c_i\ne0$ for $i=1,2,3$, Remark~\ref{REM_GeometricKEngel} ensures that there can be no geometric K-Engel structure.
On the other hand if $c_i=0$ for some $i=1,2,3$, without loss of generality we can suppose $c_1=0$, then $c_2=-c_3\ne0$ and we have 
\[
  Y:=[W,X]=c_3X_2-c_3X_3,\qquad R=X_1
\]
so that $R$ commutes with $W,\,X$ and $Y$, impliying that $\{W,\,X,\,Y,\,R\}$ is a K-Engel framing and we conclude as above.

\subsection{$Sol_0^4$}
Suppose that $P(m,n)$ given by Equation~\eqref{EQN_Pmn} has two distinct complex solutions $e^\lambda$ and $e^{\bar\lambda}$ and a real solution $e^{-2a}$, where $\Re(\lambda)=a$ and $\Im(\lambda)=b$.
The Lie algebra is generated by $U_1,\,U_2,\,V$, and $T$ satisfying
\begin{align*}
 [T,U_1]=aU_1+bU_2,\qquad[T,U_2]=-bU_1+aU_2,\qquad[T,V]=-2aV,
\end{align*}
and such that all other Lie brackets vanish.
The only geometric Engel structure here is $\sD=\bra W,\,X\ket$ where $W=U_1+V$ and $X=T$.
Notice that if a left-invariant vector field $Z=a_1U_1+a_2U_2+vV+tT$ satisfies $[Z,W]=0=[Z,X]$ then we must have
 \begin{align*}
 -a a_1 + b a_2&=0\\
 -b a_1 - a a_2&=0\\
 \ 2 a v = 2 a t &= 0
\end{align*}
but since both $a$ and $b$ must be non-zero we conclude that $Z=0$.
Remark~\ref{REM_GeometricKEngel} ensures that there can be no geometric K-Engel structure on $Sol_0^4$.

\subsection{$Sol_1^4$}
The Lie algebra is generated by $A,\,B,\,C$, and $T$ with the relations
\[
 [T,A]=-A,\qquad[T,B]=B,\qquad[A,B]=C,
\]
and such that all other Lie brackets vanish.
We take the Engel structure $\sD=\bra W,\,X\ket$ where $W=T$ and $X=A+B$, then
\[
 Y:=[W,X]=-A+B,\qquad R=[X,Y]=2C,
\]
so that $R$ commutes with $W,\,X$ and $Y$, implying that $\{W,\,X,\,Y,\,R\}$ is a K-Engel framing and we conclude as above.

\subsection{$Nil^4$}
The Lie algebra is generated by $A,\,B,\,C$, and $D$, where the only non-vanishing Lie brackets are
\[
 [D,A]=B,\qquad[D,B]=C.
\]
We take the Engel structure $\sD=\bra W,\,X\ket$ where $W=A$ and $X=D$, then
\[
 T:=[W,X]=-B,\qquad R:=[X,Y]=-C
\]
so that $R$ commutes with $W,\,X$, and $T$, impliying that $\{W,\,X,\,T,\,R\}$ is a K-Engel framing and we conclude as above.
Notice that if we take the dual coframing $\{a,\,b,\,c,\,d\}$, we can fix defining forms $\alpha=c$ and $\beta=-b$ which give Reeb distribution $\sR=\bra T,\,R\ket$.
Corollary~\ref{COR_TotallyGeoFraming} ensures that $\sD$ is totally geodesic with respect to the metric making this framing orthonormal.

\end{document}